\documentclass[english,12pt,leqno]{article}
\usepackage[latin1]{inputenc}

\usepackage{fullpage}

\usepackage{tipa}
\usepackage{dsfont}
\usepackage{tikz-cd}
\usepackage{float}
\usepackage{graphicx}

\usepackage{amsxtra}
\usepackage{amsmath}
\usepackage{amssymb}
\usepackage{amsfonts}
\usepackage[all,2cell]{xy}
\usepackage{mathrsfs}
\usepackage{amsthm}
\usepackage{enumitem}
\usepackage{hyperref}
\usepackage{subfigure}
\usepackage{mathabx}
\usepackage{euscript}
\usepackage[bbgreekl]{mathbbol}
\usepackage{pgfplots}

\usetikzlibrary{patterns}
\usetikzlibrary{calc}

\usepackage[OT2,T1]{fontenc}
\DeclareSymbolFont{cyrletters}{OT2}{wncyr}{m}{n}
\DeclareMathSymbol{\Sha}{\mathalpha}{cyrletters}{"58}

\makeatletter
\hypersetup{
unicode=true,           
pdftoolbar=true,        
pdfmenubar=true,        
pdffitwindow=false,     
pdfstartview={FitH},    
pdftitle={\@title},     
pdfauthor={\@author},   
pdfsubject={},          
pdfcreator={},          
pdfproducer={},         
pdfkeywords={},         
pdfnewwindow=true,      
colorlinks,             
linkcolor=black,        
citecolor=black,        
filecolor=black,        
urlcolor=black          
}
\makeatother

\newtheorem{thm}[equation]{Theorem}
\newtheorem{cor}[equation]{Corollary}
\newtheorem{lem}[equation]{Lemma}
\newtheorem{prop}[equation]{Proposition}
\newtheorem{refo}[equation]{Reformulation}

\newtheoremstyle{example}{\topsep}{\topsep}%
 {}
 {}
 {\bfseries}
 {.}
 {2pt}
 {\thmname{#1}\thmnumber{ #2}\thmnote{ #3}}

\theoremstyle{example}

\newtheorem{Defi}[equation]{Definition}
\newtheorem{rem}[equation]{Remark}

\newtheorem{ex}[equation]{Example}

\numberwithin{equation}{section}

\def\CC{\mathbb{C}}

\def\RR{\mathbb{R}}
\def\ZZ{\mathbb{Z}}

\def\aen{\mathfrak{a}}

\def\gen{\mathfrak{g}}
\def\hen{\mathfrak{h}}

\def\len{\mathfrak{l}}

\def\Ben{\mathfrak{B}}

\def\Pen{\mathfrak{P}}

\def\Fc{\mathcal{F}}

\def\Mc{\mathcal{M}}

\def\Hc{\mathcal{H}}

\def\Pc{\mathcal{P}}

\def\Sc{\mathcal{S}}

\def\Vc{\mathcal{V}}

\def\Alex{{\on{Alex}}}

\def\Br{{\on{Br}}}
\def\ba{{{\bf a}}}
\def\bl{{{\bf l}}}
\def\bn{{{\bf n}}}

\def\bz{{{\bf z}}}

\def\Ch{{\on{Ch}}}
\def\CM{{\on{CM}}}
\def\CMen{{\mathfrak {CM}}}

\def\del{{\partial}}

\def\eps{{\varepsilon}}

\def\Fun{{\on{Fun}}}

\def\Hom{\operatorname{Hom}\nolimits}
\def\HS{{\cal{HS}}}

\def\Id{\operatorname{Id}\nolimits}

\def\k{\mathbf k}

\def\Lex{{\on{Lex}}}
\def\lla{\longleftarrow}
\def\lra{\longrightarrow}

\def\Mor{{\on{Mor}}}

\def\Ob{\operatorname{Ob}\nolimits}
\def\on{\operatorname}
\def\ol{\overline}
\def\oo{{\infty}}
\def\OP{{\mathcal{OP}}}

\def\Perv{{\on{Perv}}}
\def\phi{{\varphi}}
\def\Ps{{\on{Ps}}}

\def\Ran{{\on{Ran}}}

\def\Sn{{\on{Sn}}}
\def\Sup{{\on{Sup}}}

\def\Sym{{\on{Sym}}}

\newcommand{\ta}{\tilde a}

\def\ul{\underline}

\def\Vect{\on{Vect}}

\def\wh{ \widehat}
\def\wt{\widetilde}

\def\+{{\oplus}}
\def\- {{\setminus}}
\def\x{{\otimes}}
\def\1{{\mathbf{1}}}
\def\2{{\mathbf{2}}}
\def\(({(\hskip -1mm (}
\def\)){)\hskip -1mm )}
\def\be{\begin{equation}}
\def\ee{\end{equation}}
\def\ed{\end{document}}

\title{  PROBs  and perverse sheaves I. Symmetric products}
\author{  Mikhail Kapranov, Vadim Schechtman}

\begin{document}

\maketitle


 \begin{abstract}
 Algebraic structures involving both multiplications and comultiplications (such as, e.g., bialgebras or Hopf
 algebras) can be encoded using PROPs (categories with PROducts and Permutations) of Adams and MacLane.
 To encode such structures on objects of a braided monoidal category,  we need PROBs (braided analogs of
 PROPs). Colored PROBs correspond to multi-sorted structures. 
 
 In particular, we have a colored PROB $\Ben$ governing $\ZZ_{\geq 0}$-graded bialgebras 
 in braided categories. As a category,
 $\Ben$ splits into blocks $\Ben_n$ according to the grading. We relate $\Ben_n$ with the category $\Pen_n$
 of perverse sheaves on the symmetric product $\Sym^n(\CC)$ smooth with respect to the natural stratification
 by multiplicities. More precisely, we show that $\Pen_n$ is equivalent to the category of functors $\Ben_n\to\Vect$.
 This gives a natural quiver description of $\Pen_n$. 
  \end{abstract}



\section{Introduction. The main result}\label{sec:main-res}

Let $\k$ be a field. All vector spaces  in this paper  will be assumed
 $\k$-vector spaces  and 
 all  categories  will be assumed
 $\k$-linear. 

\paragraph{PROPs and PROBs.} The concept of a PROP  was introduced by Adams and MacLane \cite{adams, maclane}. It allows one
 to axiomatize algebraic structures on a vector space (or, more generally, on
an object $A$ of a symmetric monoidal category) involving both multiplications
$A^{\otimes m} \to A$ and  comultiplications $A\to A^{\otimes n}$. 
See \cite{markl-prop} for a modern
exposition. 

\vskip .2cm

The term PROP is an abbreviation for  ``category with PROducts and Permutations''.
 Explicitly, a PROP is a symmetric monoidal category
  $(\Pc, \otimes, \1)$ with $\Ob(\Pc) = \{ [m], m\in\ZZ_+\}$ identified with the set of non-negative integers so that the tensor operation
$\otimes$ is, on objects, given by the addition of integers: $[m]\otimes [n]=[m+n]$
and the unit object is  $\1=[0]$.  In other words, a PROP is a strict symmetric monoidal
category  whose objects are tensor powers of a single object, denoted $[1]$.

Given a symmetric monoidal category $\Vc$, we can speak
about 
 {\em algebras}   in $\Vc$ over a PROP $\Pc$. Such an algebra
 is simply a symmetric monoidal functor $F: \Pc \to \Vc$. If we denote $A=F([1]) \in\Vc$,
 then each space $\Hom_\Pc([m], [n])$ is mapped into the space of mixed operations
 $A^{\otimes m} \to A^{\otimes n}$. 
 
 \begin{ex}[(The PROP of Hopf algebras)] A basic example of a structure involving multiplications and comultiplications is that of a Hopf algebra. A Hopf algebra in
 a symmetric monoidal category $\Vc$ is an object $A$ together with a multiplication
 $\mu: A\otimes A\to A$, comultiplication $c: A\to A\otimes A$, unit $e: \1\to A$ and counit
 $\eta: A\to\1$ satisfying the well known relations. The corresponding PROP, denote it
 $\HS$, can be seen as the symmetric monoidal category generated by the {\em universal
 Hopf algebra} $\ba=[1]$. 
 That is,  for any symmetric monoidal category $\Vc$,  we have a bijection between:
\begin{itemize}
\item[(i)] Hopf algebras $A$ in $\Vc$;

\item[(ii)]  Symmetric monoidal functors $F: \HS\to\Vc$,
\end{itemize}   
given by $A= F(\ba)$. At a more intuitive level, 
  $\HS$ contains the morphisms
 $\mu: \ba\otimes\ba\to\ba$ and similarly $c,e,\eta$ as above, together with all their iterated compositions and tensor 
 products, which are subject to the relations of a Hopf algebra ``and nothing else''. 
 
Along with Hopf algebras, we can consider more general objects, namely {\em bialgebras}
in symmetric monoidal categories. Thus a bialgebra $A$ has compatible associative $\mu: A\otimes A\to A$ and coassociative $c: A\to A\otimes A$ but the unit and counit are not required. As before, we have the PROP $\HS^+$ describing bialgebras. It was introduced by Markl \cite{markl}. 
  \end{ex}

Further, it is well known
\cite{majid, takeuchi-hopf} the concept of a Hopf algebra  (or, more generally of a bialgebra)
can
be defined in any braided, not necessarily symmetric,  monoidal category $\Vc$. 
To study such structures we need to modify the concept of a PROP to that of
a PROB (with Braidings instead of Permutations). Thus, a PROB is a braided
monoidal category whose objects are tensor powers of a single object $\ba=[1]$. 

\begin{ex}[(The PROB of braided Hopf algebras)]
As before, there is a PROB $\Hc$ such that Hopf algebras in a braided monoidal category
$\Vc$ are in bijection with braided monoidal functors $\Hc\to\Vc$. Its objects are
tensor powers of the {\em universal braided Hopf algebra} $\ba=[1]\in\Hc$ and morphisms are
iterated compositions and tensor products of the generating morphisms $\mu,c,e,\eta$
as above subject to the relations of a braided Hopf algebra. This PROB was 
introduced by Habiro  \cite{habiro-bottom}  who denoted it  $\langle {\tt H}\rangle$. 
\end{ex}

Despite their deceptively short definitions by universal properties, the PROP $\HS$
and the PROB $\Hc$ are quite non-trivial objects. In particular, their stucture as
ordinary categories, i.e., some description of the spaces $\Hom([m], [n])$
for all $m,n$, is not so easy to pin down. For the simpler PROP $\HS^+$
  such a description was obtained  by Pirashvili \cite{pirashvili}. 

\vskip .2cm

The goal of this paper and the one to follow \cite{KS-Ran} is to relate  versions of the PROB $\Hc$
to  quivers describing perverse sheaves on certain configuration spaces. Let us
describe the version relevant for this paper.

\paragraph{The PROP of graded bialgebras.}\label{par:uni-gr-bi}
By a {\em graded bialgebra} in a braided  monoidal category $\Vc$
we mean a bialgebra decomposed into a direct sum $A =\bigoplus_{n\geq 0} A_n$  so that 
$A_0=\1$ is the unit object, the multiplication and comultiplication are homogeneous and their
components involving $A_0$ are the identities. See \cite{KS-shuffle} where such objects were
called primitive bialgebras.  A graded bialgebra automatically has a unit, counit
and   antipode, see \cite{KS-shuffle}, Prop. 2.4.11. 

\vskip .2cm

In order to be able to speak about graded bialgebras, the existence of direct sums in 
$\Vc$ is, strictly speaking, not necessary as all the conditions can be
reformulated in terms of the individual 
graded components $A_n$. Therefore we will
understand a graded bialgebra $A$ as a collection of these components:
$A=(A_n)_{n\geq 0}$. The direct sum 
$\bigoplus_{n\geq 0} A_n$  can be always considered as an object
of the formal direct sum completion of $\Vc$ but we need not require that it
belongs to $\Vc$.

\vskip .2cm

As before,
there is a braided category $\Ben$  generated by
 the {\em universal graded bialgebra}
$\ba= (\ba_n)_ {n\geq 0}$, $\quad \ba_0=\1$. 
Objects of $\Ben$ are formal tensor products
$\ba_\alpha = \ba_{\alpha_1}\otimes \cdots\otimes\ba_{\alpha_r}$
associated to all the {\em ordered partitions}, i.e., sequences of positive integers 
$\alpha = (\alpha_1, \cdots, \alpha_r)$. 
Morphisms are
generated by the  elementary formal morphisms
\[
\mu_{p,q}:\ba_p\otimes \ba_q\to  \ba_{p+q}, \quad \Delta_{p,q}: \ba_{p+q}\to 
\ba_p\otimes \ba_q,\quad \mu_{p,0}=\mu_{0,p}=\Delta_{p,0}=\Delta_{0,p}=\Id_{\ba_p},
\]
as well as the braidings, modulo the relations following from the axioms of
a graded bialgebra and a braided category. 

The braided category $\Ben$ is an example of a {\em colored PROB} in that it describes
algebraic structures not on a single object but on a family of objects $(A_n)$ of a braided category. 
 For a discussion of colored PROPs, not PROBs see  \cite{hackney-robertson}.

\vskip .2cm

Let $\OP$ be the set of all  ordered partitions $\alpha$ as above and $\Pc\subset \OP$ be the set of
{\em unordered} (or classical) {\em partitions}, i.e., sequences $\alpha = (\alpha_1\geq \cdots\geq \alpha_r > 0)$.
The $\ba_\alpha$, $\alpha\in\Pc$ form a system of representatives of isomorphism classes of objects of $\Ben$. 

\vskip .2cm

For $\alpha\in\OP$
let $|\alpha|=\sum \alpha_i$. 
Let $\OP_n\subset\OP$, resp. $\Pc_n\subset\Pc$  consist of $\alpha$
such that $|\alpha|=n$, i.e., that $\alpha$ is an ordered resp. unordered partition of $n$. 
Let $\Ben_n\subset\Ben$ be the full subcategory on objects $\ba_\alpha$
with $\alpha\in\Pc_n$. 
It is not closed under the product $\otimes$. The full subcategory on the $\ba_\alpha$, $\alpha\in\OP_n$, is equivalent to $\Ben_n$.

\paragraph{Symmetric products and perverse sheaves.} Let $\Sym^n(\CC)$
be the $n$th symmetric product of $\CC$, i.e., the space of monic
polynomials
\[
f(x) \,=\, x^n + a_1 x^{n-1} +\cdots + a_n,\quad a_i\in\CC
\]
or, equivalently, the space of effective divisors $\bz = \sum_i n_i z_i$ with 
$z_i\in\CC$, $n_i\geq 0$, $\sum n_i=n$. Each such divisor $\bz$ has a {\em type}
which is the partition
$ t(\bz)\in\Pc_n$ obtained by 
arranging the $n_i$ in a non-increasing
order, and we denote  $\Sym^\alpha(\CC)\subset \Sym^n(\CC)$
the subspace of $\bz$ of type $\alpha$. This gives an algebraic
Whitney stratification of $\Sym^n(\CC)$ which we denote $\Sc^{(0)}$
and call the {\em stratification by multiplicities}. The open strarum of $\Sc ^{(0)}$ is 
\[
\Sym^n_\neq(\CC) = \Sym^{(1,\cdots, 1)}(\CC), 
\]
the space of multiplicity-free divisors or, equivalently, of polynomials $f(x)$ with non-zero discriminant. 

\vskip .2cm

Let $\Vc$ be any abelian category (not assumed   monoidal).
We can then speak about $\Vc$-valued  perverse sheaves
(with respect to the middle perversity) on $\Sym^n(\CC)$
which are constructible with respect to the stratification $\Sc^{(0)}$,
see \cite{KS-shuffle}. They form an abelian category   which
we denote $\Perv(\Sym^n(\CC); \Vc)$. 
For example, if $\Vc=\Vect_\k$ is the category of $\k$-vector
spaces, then  $\Perv(\Sym^n(\CC); \Vc)$ is the category of perverse
sheaves of $\k$-vector spaces in the usual sense. 
Here is our  main result, whose  proof   will be given at the end of  Section \ref{sect:CM-Ben}.

\begin{thm}\label{thm:main-sym}
 We have an equivalence of categories $\Perv(\Sym^n(\CC);\Vc)\simeq
\Fun(\Ben_n, \Vc)$. 

\vskip .2cm

\end{thm}

In other words, we have  an elementary, or   quiver  description
of the category of  perverse sheaves on $\Sym^n(\CC)$.  The corresponding quiver (with relations) is
the category $\Ben_n$. The vertices of this quiver are   the  objects of $\Ben_n$, 
i.e., the
$\ba_\alpha$, $\alpha\in\Pc_n$. They are in bijection with the strata $\Sym^\alpha(\CC)$
of the stratification $\Sc^{(0)}$. 

\begin{rem}
 The monoidal  structure $\otimes_{m,n}: \Ben_m\times\Ben_m\to\Ben_{m+n}$ corresponds, 
 at the level of perverse sheaves, to the  functor (``comonoidal structure'')  
\[
\nabla_{m,n}: 
\Perv(\Sym^{m+n}(\CC); \Vc) \lra \Perv(\Sym^m(\CC)\times\Sym^n(\CC); \Vc)
\]
defined  geometrically as  follows. Choose two disjoint open disks $D, D'\subset \CC$ so that we have
open embeddings
\[
\Sym^m(\CC) \times\Sym^n(\CC) \buildrel i\over\lla \Sym^m(D) \times\Sym^n(D')
\buildrel j\over\lra \Sym^{m+n}(\CC)
\]
such that the pullback functor $i^*$ defines an equivalences on the categories of
perverse sheaves with respect to the natural stratifications. Then 
$\nabla_{m,n}= (i^*)^{-1} \circ j^*$. 
 
\end{rem}

 \paragraph{The simplest example.}
(1) Let $n=2$. The category $\Ben_2$ has two objects, $\ba_1\otimes \ba_1$ and $\ba_2$, and its
morphisms are generated by: 
\[
\xymatrix{
\ar@(ul,dl)_R\hskip -1cm& \hskip -1cm \ba_1\otimes\ba_1 \ar@<.7ex>[r]^{ \mu_{1,1}}&\ba_2
\ar@<.7ex>[l]^{\hskip .2cm \Delta_{1,1}}   
}
\]
with $R$ being the braiding, therefore invertible. Note that $\ba_1$ is primitive, i.e.,
\[
\Delta|_{\ba_1} \,=\,\Delta_{1,0} + \Delta_{0,1}: \ba_1 \lra (\ba_1\otimes\1)\, \oplus \, (1\otimes\ba_1)
\]
is the sum of two copies of the identity morphism of $\ba_1$. This implies that $\Delta_{1,1}\mu_{1,1} = \Id + R$, 
cf.  \cite{KS-shuffle}, \S 5.2. In particular, $\Id-\Delta_{1,1}\mu_{1,1} = -R$ is invertible. This means that
$\Fun(\Ben_2,\Vc)$ is the category of diagrams in $\Vc$ 
\[
\xymatrix{
 \Psi
\ar@<.7ex>[r]^{b}&\Phi
\ar@<.7ex>[l]^{a}  
}
\]
such that $T_\Psi=\Id_\Psi-ab$ is invertible (which implies that $T_\Phi=\Id_\Phi-ba$ is invertible, see Eq. (1.1.6) of
\cite{KS-shuffle}). 

\vskip .2cm

(2) On the other hand, $\Sym^2(\CC)=\CC^2$ is the space of quadratic polynomails $x^2+a_1x+a_2$, and the
only non-open stratum of $\Sc^{(0)}$ is the parabola given by $a_1^2=4a_2$ formed by polynomials with a double
root. Factoring out by translational symmetry, we find that
\[
\Perv(\Sym^2(\CC);\Vc) \,\simeq \, \Perv(\CC,0; \Vc)
\]
is identified with the category of $\Vc$-valued perverse sheaves on $\CC$ with the only possible singularity at $0$.

\vskip .2cm

The $n=2$ instance of Theorem \ref{thm:main-sym} reduces therefore to (a $\Vc$-valued version of) the classical result
of \cite{beil-gluing, {galligo-GM}} describing $\Perv(\CC,0;\Vect_\k)$ in terms of $(\Phi, \Psi)$-diagrams
of vector spaces as above.

  \paragraph{ Discussion and further plans. }  Theorem \ref{thm:main-sym} can be seen as a
  refinement of two previous results:
  
  \vskip .2cm
  
  (1)  The main result of \cite{KS-shuffle} which identifies the category of graded bialgebras  $A=(A_n)_{n\geq 0}$ 
  in  a braided monoidal $\Vc$  
 with  that of {\em factorizable systems} $(\Fc_n)_{n\geq 0}$ of perverse sheaves
   on all the symmetric products
  $\Sym^n(\CC)$. 
   The refinement consists in passing from such factorizable systems  to individual
  perverse sheaves on an individual symmetric product and in allowing $\Vc$ to be 
  an arbitrary abelian (not necessarily monoidal) category. 
  
  \vskip .2cm
  
  (2) The special case $\gen = \gen\len_n$ of the main result of  \cite{KS-hW} 
 which describes perverse sheaves on any quotient
  $W\backslash \hen$ where $\hen$ is the Cartan subalgebra of a complex reductive Lie algebra $\gen$
  and $W$ is the Weyl group. If $\gen=\gen\len_n$, then $W\backslash \hen=\Sym^n(\CC)$,
  and we get a description of $\Perv(\Sym^n(\CC))$. This description, however,  is more cumbersome than
   Theorem \ref{thm:main-sym} so the present refinement consists in giving a neater one-shot description. 
   
   \vskip .2cm
  
Our proof of   Theorem \ref{thm:main-sym}   uses the results (1) and (2) above 
  by assembling the categories constructed in (2) into a single braided category $\CMen$
  and constructing a graded bialgebra in  this category so that application of (1)  leads to an
  identification $\CMen\sim\Ben$. 
  
  It is also interesting to note the similarity between the description of
  the PROP $\HS^+$ given by Pirashvili \cite{pirashvili} and the description (2) above proceeding
  in terms of so-called {\em mixed Bruhat sheaves} \cite{KS-hW}.  Both descriptions involve natural ``bivariant''
  objects: covariant in one direction, contravariant in the other with some  base change-type relations
  relating the two variances. 
  
  \vskip .2cm
  
  In a  sequel to this paper  \cite{KS-Ran} we plan to describe the category of perverse sheaves on the 
  {\em Ran space}  $\Ran(\CC)$ in terms of a category related to the PROB $\Hc$ governing
  braided Hopf algebras. 
  
  \vskip .2cm
  
  Both papers can be seen as developing the observation, going back to Lurie, that bialgebras
  are Koszul dual to $E_2$-algebras and thus \cite{CG, lurie-HA} to locally constant factorization algebras
  on $\RR^2=\CC$ i.e., to factorizable (complexes of) sheaves on the Ran space.
   Informally, an $E_2$-algebra can be seen as a cochain complex $E$
  with two (homotopy) compatible (homotopy) associative multiplications. Now, Koszul duality
 gives an equivalence between associative dg-algebras and coalgebras. Applying it
 to one of the two multiplications on $E$, we get a structure consisting  of  (homotopy)
 compatible multiplication and
 comultiplication, i.e., a homotopy version of a dg-bialgebra $E^!$. Koszul duality being a derived equivalence, for $E^!$ to be an
 honest (non-dg) bialgebra,  $E$ must be a nontrivial complex. What makes our approach work is
 a remarkable match between this type of complexes and 
 the Cousin complexes playing an essential role in our earlier descriptions
 of perverse sheaves \cite{KS-shuffle, KS-hW}. 
 
 \paragraph{Outline of the paper.} Apart from the present introductory \S \ref {sec:main-res}, 
 the paper has three more sections. 
 
 \vskip .2cm
 
 In \S \ref{sec:PROB-cont} we, first,  
  give  background material on braided categories and graded bialgebras in such categories. 
  In particular, we give a self-contained treatment of Deligne's interpretation
  of braidings in terms of ``$2$-dimensional tensor products'' of objects labelled by 
  points in the plane. 
 We recall the concept of contingency matrices and their vertical and horizontal contractions
 from \cite{KS-cont}. We further  associate to
  a contingency matrix $M=\|m_{ij}\|$  and a graded bialgebra $A$ an object $A_M$  which is the $2$-dimensional tensor product of the components
  $A_{m_{ij}}$. We use the multiplication and comultiplication in $A$ to connect the objects $A_M$, and $A_N$ whenever $M$ is obtained from $N$   by a vertical or horizontal contractions and establish (Proposition \ref{prop:gr-bi-Hodge}) a system of relations
  for such connecting morphisms. 
  
  \vskip .2cm
  
  In \S \ref{sec:cont-mat-cat} we define a category $\CMen$  
  whose objects are symbols $[M]$ associated to contingency matrices $M$ and 
  the relations of Proposition \ref{prop:gr-bi-Hodge} are promoted into
   into a system of defining
  relations for  the morphisms of the category.  Thus any graded bialgebra $A$ in any braided category
   $\Vc$ gives rise to a functor
  $\xi_A: \CMen\to\Vc$ (Corollary \ref{cor:bialg-CM-V-1}) sending $[M]$ to $A_M$. 
  We further make $\CMen$ into a braided monoidal category so that
 the functor $\xi_A$ above is in fact braided monoidal (Proposition \ref{prop:bialg-CM-V}). 
 
 \vskip .2cm
 
 In \S \ref{sect:CM-Ben} we notice that $\CMen$ carries a braided bialgebra $\aen$
 with components $[n]$ associated to $1\times 1$ matrices $n$. This allows us to
 connect $\Ben$ and $\CMen$ by a braided monoidal functor which we show to be an
 eqivalence (Theorem \ref{thm:CM=B}). 
 Finally, for each $n$ we compare the degree $n$ block $\CMen_n$
 with the specialization, for $\gen=\gen\len_n$, of the concept of mixed Bruhat sheaf
 which was introduced in  \cite{KS-hW}  for description of perverse sheaves
 on the adjont quotient $W\backslash \hen$ of a reductive Lie algebra $\gen$. 
 This comparison yields the first in the two identifications below
 \[
 \Perv(\Sym^n(\CC);\Vc) \,  \simeq \, \Fun(\CMen_n, \Vc) \,
 \buildrel \text{Thm. \ref{thm:CM=B}}\over
 \simeq \, 
 \Fun(\Ben_n, \Vc)
 \]
 thus proving Theorem \ref{thm:main-sym}.

 \paragraph{Acknowledgements.} The research of M.K. was supported by the  World Premier International Research Center Initiative (WPI Initiative), 
 MEXT, Japan.

\vfill\eject

\section {Graded bialgebras   and contingency matrices}\label {sec:PROB-cont}

\paragraph{ Braids, braided categories and bialgebras.}
 
Let $\Br_n$ be the braid group on $n$ strands, with the standard generators
$\sigma_1,\cdots, \sigma_{n-1}$ and relations
\be\label{eq:br-relations}
\sigma_p\sigma_{p+1}\sigma_p = \sigma_{p+1}\sigma_p \sigma_{p+1},
\quad \sigma_p \sigma_q=\sigma_q\sigma_p, \,\,\,|p-q|\geq 2.
\ee
Let also $S_n$ be the symmetric group with the standard
generators $\ol\sigma_1,\cdots,\ol\sigma_{n-1}$ subject to the same relations
as in \eqref{eq:br-relations}
 together with $\ol\sigma_i^2=1$. Thus we have the surjective morphism
 $p: \Br_n\to S_n$.

  \vskip .2cm

 By a {\em  monoidal category} we will  mean a strictly associative
 monoidal category, with 
a strict unit object denoted by $\1$. Let $\Vc$ be a braided monoidal category with braiding denoted by $R=(R_{V,W}: V\otimes W\to W\otimes V)$.

 \vskip .2cm
Given objects
$V_1,\cdots, V_n$ of $\Vc$, and an element $b\in \Br_n$, we have
the permutation $t=p(b)\in S_n$, and the braiding isomorphism
\[
R_b: V_1\otimes\cdots\otimes V_n\lra V_{t(1)}\otimes\cdots\otimes
V_{t(n)}. 
\]

A {\em bialgebra} in $\Vc$  is an
object $A$  of $\Vc$ equipped with an associative multiplication $\mu: A\otimes A\to A$ and a
coassociative comultiplication $\Delta: A\to A\otimes A$ satisfying the following
compatibility condition: $\Delta$ is a morphism of algebras, if the multiplication on $A\otimes A$ is defined using the braiding:
\[
A\otimes A \otimes A\otimes A \buildrel \Id\otimes R_{A,A}\otimes \Id \over\lra  A\otimes A \otimes A\otimes A 
\buildrel \mu\otimes\mu \over\lra A\otimes A. 
\] 

\paragraph {Geometric interpretation of braided categories.} 
Let $\Vc$ be a category, $I$  be a finite set and $(V_i)_{i\in I}$
be a family of objects
of $\Vc$ labelled by $I$. If $\Vc$ is symmetric monoidal, then we can speak
about the object $\bigotimes_{i\in I} V_i$ without specifying an order on $I$:
the ordered tensor products in different orders are canonically identified with
each other.

\vskip .2cm

If $\Vc$ is braided monoidal, then  the notation 
$\bigotimes_{i\in I} V_i$ does not make
sense, as there is no single canonical identification between two given
ordered products. It was pointed out by Deligne that the correct
structure on $I$ to make the product canonical is not an ordering
(which of course suffices) but an embedding $I\hookrightarrow  \CC$. 
In other words, once we assign to each $i\in I$ a complex number $z_i\in\CC$
such that $z_i\neq z_j$ for $i\neq j$, there is a ``well-defined object''
$\bigotimes_{i\in I} V_i(z_i)$ (the $2$-dimennsional tensor product with respect to $V_i$
positioned as $z_i$). When the $z_i$ move, these objects unite into a local
system on $\CC^I_\neq$ whose monodromy gives the braiding.
For convenience of the reader we recall a precise elementary construction. 

\begin{Defi}
Let $\Vc$ be a category.

\vskip .2cm

(a) A {\em pseudo-object} (or an object defined up to a unique isomorphism)
of $\Vc$ is a datum $V$ of a set $K$, of objects $V_k\in\Vc$ for each $k\in K$
and of isomorphisms $\phi_{k,k'}: V_k\to V_{k'}$ given for each $k,k'$
and satisfying
\[
\phi_{kk}=\Id,\quad \phi_{k',k''} \phi_{k,k'} = \phi_{k, k''},\quad \forall\,\,k, k', k''. 
\]
The objects $V_k$ are called the {\em determinations} of the pseudo-object $V$,
and the morphisms $\phi_{kk'}$ are called the {\em transition maps} of $V$. 

\vskip .2cm

(b) A morphism $u$ from a pseudo-object $V=(V_k, \phi_{k,k'})_{k,k'\in K}$
to a pseudo-object $W=(W_l, \psi_{l,l'})_{l,l'\in L}$ of $\Vc$ is a datum of morphisms
$u_{k,l}: V_k\to W_l$ for all $k\in K$, $l\in L$ such that
\[
u_{k,l}=u_{k',l} \phi_{k,k'}, \quad u_{k,l'}=\psi_{l,l'}u_{k,l},
\quad \forall \,\, k,k'\in K,\,\, l,l'\in L. 
\]
\end{Defi}

With this definition, pseudo-objects in $\Vc$ form a category $\Ps(\Vc)$.
Any actual object of $\Vc$ can be considered as a pseudo-object with $|K|=1$.
This gives a functor $\Vc\to\Ps(\Vc)$ which is easily seen to be an equivalence.
In this way a pseudo-object can be seen to be  ``as good as an actual object''
of $\Vc$. 

\vskip .2cm

We will construct $\bigotimes_{i\in I} V_i(z_i)$ as a pseudo-object and
start with describing its indexing set $K$. 

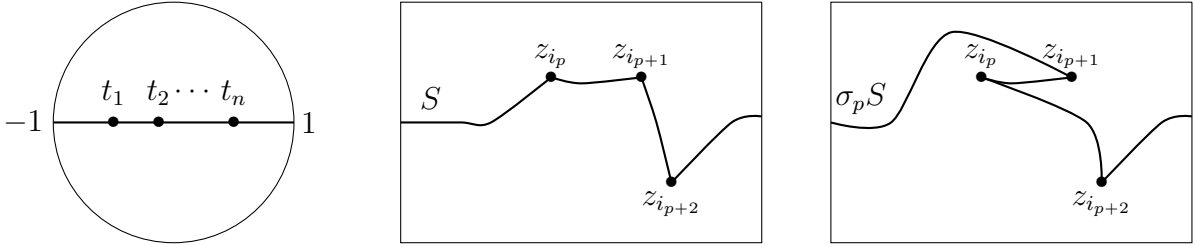
\begin{figure}[H]
\centering
\begin{tikzpicture}[scale=0.4]
\draw (0,0) circle (4); 
\draw [line width=0.8] (-4,0) -- (4,0); 
 \node at (-2,0){\small$\bullet$}; 
   \node at (-0.5,0){\small$\bullet$}; 
     \node at (2,0){\small$\bullet$}; 
\node at (-2,1){$t_1$};  
\node at (-0.5, 1){$t_2$};         
\node at (2,1){$t_n$};  
\node at (0.7, 1){$\cdots$};   
\node at (-5,0){$-1$}; 
\node at (4.5,0){$1$}; 

\end{tikzpicture}
\quad\quad 
\begin{tikzpicture}[scale=0.4]

\draw (6,4) --(-6,4) -- (-6,-4) -- (6,-4) -- (6,4);   
\node at (-1,1.5){\small$\bullet$};
\node at (2,1.5){\small$\bullet$};
\node at (3,-2){\small$\bullet$}; 

\draw [line width=0.8] plot [smooth, tension=0.5]
coordinates  {(-6,0) (-4,0) (-3,0) (-1, 1.5) }; 

\draw  [line width=0.8]  plot [smooth, tension=0.5]
coordinates { (-1,1.5) (0,1.3) (2,1.5)} ; 

\draw   [line width=0.8]  plot [smooth, tension=0.5]
coordinates {(2,1.5) (2.5, 0) (3,-2)} ; 

  \draw  [line width=0.8]   plot [smooth, tension=0.5]
coordinates {  (3,-2) (5,0) (6,0.2) } ; 

\node at (-5,0.8){$S$}; 
\node at (-1,2.2){$z_{i_p}$}; 
\node at (2,2.2){$z_{i_{p+1}}$};  
\node at (3, -2.7){$z_{i_{p+2}}$};

\end{tikzpicture}\quad\quad
\begin{tikzpicture}[scale=0.4]

\draw (6,4) --(-6,4) -- (-6,-4) -- (6,-4) -- (6,4);   
\node at (-1,1.5){\small$\bullet$};
\node at (2,1.5){\small$\bullet$};
\node at (3,-2){\small$\bullet$}; 

\draw [line width=0.8] plot [smooth, tension=0.5]
coordinates  {(-6,0) (-4,0) (-2,3) (2, 1.5) }; 

\draw  [line width=0.8]  plot [smooth, tension=0.5]
coordinates { (-1,1.5) (0,1.3) (2,1.5)} ; 

\draw   [line width=0.8]  plot [smooth, tension=0.5]
coordinates {(-1,1.5) (2.5, 0) (3,-2)} ; 

 \draw  [line width=0.8]   plot [smooth, tension=0.5]
coordinates {  (3,-2) (5,0) (6,0.2) } ; 

\node at (-5,0.8){$\sigma_p S$}; 
\node at (-1,2.2){$z_{i_p}$}; 
\node at (2,2.2){$z_{i_{p+1}}$};  
\node at (3, -2.7){$z_{i_{p+2}}$};

\end{tikzpicture}

\caption{Snakes as the images of the interval $(-1,1)$ in the disk.}\label{fig:snake}
\end{figure}

\begin{Defi}
Let $Z\subset \CC$ be a finite subset, $|Z|=n>0$. A {\em $Z$-snake}
is a simple curve $S\subset \CC$ which is a finite perturbation of $\RR$
and passes through each element of $Z$ once. See the center and right
of Fig. \ref{fig:snake}. 
\end{Defi}

Denote by $\Sn(Z)$ the set of isotopy classes of $Z$-snakes. 
A snake being oriented from $(-\oo)$ to $(+\oo)$, each $S\in \Sn(Z)$
gives an ordering $i_1=i_1(S), \cdots, i_n=i_n(S)$ of the set $I$. 
It  is classical that $\Sn(Z)$ is a left  torsor over $\Br_n$, with $\sigma_p S$ being
obtained from $S$ by an ``upper twist'' reversing the path between
$z_{i_p(S)}$ and $z_{i_{p+1}(S)}$, see the right of  Fig. \ref{fig:snake}. 
Conceptually, this follows from the interpretation of $\Br_n$ as the
mapping class group of the $n$-pointed disk \cite{ dehornoy, farb}. 
More precisely, let $D=\{|z|\leq 1\}$ be the closed unit disk in $\CC$ with $\del D=S^1$ the unit circle and let $-1 <t_1 <\cdots < t_n<1$
be real numbers. Then $\Br_n \simeq \pi_0(\Hc)$,  where
\[
\Hc = \on{Homeo}(D, \Id_{\del D}, \{t_1,\cdots,
t_n\})
\]
is the group of homeomorphisms of $D$, identical on $\del D$ and preserving
$\{t_1,\cdots, t_n\}$ as a set. 

\vskip .2cm

At the same time, let $\ol\CC=\CC\cup S^1_\oo$ be the disk compactification
of $\CC$ by $S^1_\oo$, the circle of directions at $\oo$. 

\begin{prop}
We have an identification $\Sn(Z)\simeq \pi_0(\Mc)$, where
$\Mc$ is the space of homeomorphisms $f: D\to\ol \CC$ which restrict to
the standard identification $\del D=S^1\to S^1_\oo$ on the boundary and
take $\{t_1,\cdots, t_n\}$ to $Z$. The snake corresponding to $f\in\Mc$
is the image $f((-1,1))$ (flattened near $\pm\oo$ to coincide with $\RR$ there). 
See Fig. \ref{fig:snake}. 
\end{prop}

\noindent{\sl Proof:} This is equivalent to the classical encoding of braids by curve
diagrams, see \cite{dehornoy} \S 1.3.3, esp. Fig. 2 there. \qed

\vskip .2cm

As $\Mc$ is a torsor over $\Hc$, we have that $\Sn(Z)=\pi_0(\Mc)$ is
a torsor over $\Br_n=\pi_0(\Hc)$.

 \begin{Defi}\label{def:V-i-z-i}

Let $\Vc$ be a braided category, $(V_i)_{i\in I}$, $|I|=n$ be a family of objects
and $z_i\in\CC$, $i\in I$ be distinct complex numbers. Let
$Z=\{z_i\}_{i\in I} \subset \CC$. We define a pseudo-object $\bigotimes_{i\in I}
V_i(z_i)$ of $\Vc$ with the indexing set $K=\Sn(Z)$. For  
$S\in \Sn(Z)$ the corresponding determination   is defined as
\[
\biggl(\bigotimes_{i\in I} V_i(z_i)\biggr)_S \,=\, V_{i_1(S)}\otimes\cdots\otimes
V_{i_n(S)},
\]
the ordered tensor product along the snake. Given two snakes $S,S'\in\Sn(Z)$
with $S'=b S$, $b \in\Br_n$, we define the transition map
\[
\phi_{S,S'}:= R_b: V_{i_1(S)}\otimes\cdots\otimes
V_{i_n(S)} \lra V_{i_1(S')}\otimes\cdots\otimes
V_{i_n(S')}
\]
to be the braiding isomorphism associated to $\beta$. 
\end{Defi}

 \paragraph{Factorization algebra point of view on a braided category.} It is
 convenient to extend Definition \ref{def:V-i-z-i} slightly, bringing it close
 to the formalism of factorization algebras \cite{CG}. 
 
 \vskip .2cm
 
 By a {\em closed disk} in $\CC$ we mean a subset $D\subset \CC$ 
which is either homeomorphic to the standard disk $\{ |z| \leq 1\} \subset \CC$ 
or
is a single point (a disk of radius $0$). 
Let $(V_i)_{i\in I}$ be as above and $Z_i\subset \CC$ be disjoint closed 
disks. We define
\[
\bigotimes_{i\in I} V_i(Z_i) \, := \, \bigotimes_{i\in I} V_i(z_i),\quad 
\forall\, z_i\in Z_i, \, i\in I. 
\]
Formally, we  view the objects in the 
RHS of this definition as the stalks of a $\Vc$-valued local system 
on the contractible space  $\prod_{i\in I} Z_i$ and we define the LHS
as the object of global sections of this local system. 

\vskip .2cm

Alternatively, let $Z=\bigcup_{i\in I} Z_i$. Then we can speak about
$Z$-{\em snakes} that is, simple curves $S\subset \CC$ which are
finite perturbations of $\RR$ which intersect each $Z_i$ along a
closed interval (possibly reducing to a single point). A  $Z$-snake
defines an ordering on $I=\pi_0(Z)$. 
 As before, 
the set $\Sn(Z)$ of isotopy classes of $Z$-snakes  is a $\Br_n$-torsor, and
 we can define $ \bigotimes_{i\in I} V_i(Z_i)$ as a pseudo-object
 with indexing set $\Sn(Z)$ consisting of ordered tensor products
 along the snakes. 
 
 \begin{prop}\label{prop:braid-ass}
 (1) Let $\Vc$ be a braided category, $f: J\to I$ a surjection of finite sets
 and $(V_j)_{j\in J}$ a family of objects of $\Vc$.  Let $(Z_i)_{i\in I}$
 and $(Y_j)_{j\in J}$ be two families of closed disks in $\CC$,
 each consisting of disjoint disks and such that $Y_j\subset Z_{f(j)}$, $j\in J$.
 In each such situation we have a canonical {\em associativity isomorphism}
 \[
 \alpha_f: \bigotimes_{i\in I} \biggl( \bigotimes_{j\in f^{-1}(i)} V_j(Y_j)\biggr)(Z_i) 
 \lra \bigotimes_{j\in J} V_j(Y_j). 
 \]
 These isomorphisms satisfy the following compatibility.
 
 \vskip .2cm
 
 (2) Let $K\buildrel g\over\to J\buildrel f\over \to I$ be two composable
 surjections of finite sets, and $(V_k)_{k\in K}$ be a family of objects of $\Vc$.
 Let $(X_k)_{k\in K}$, $(Y_j)_{j\in J}$ and $(Z_i)_{i\in I}$ be three
 families of closed disks in $\CC$, each consisting of disjoint disks and such
 that $X_k\subset Y_{g(k)}$, $k\in K$ and $Y_j\subset Z_{f(j)}$, $j\in J$.
 For each $i\in I$ let $g_i: (fg)^{-1}(i)\to f^{-1}(i)$ be the restriction of $g$. 
 Then the following diagram is commutative:
 \[
 \xymatrix{
 \bigotimes\limits_{i\in I} \biggl( \bigotimes\limits_{j\in f^{-1}(i)}\biggl( \bigotimes\limits_{k\in g^{-1}(j)}
 V_k(X_k)\biggr)(Y_j)\biggr)(Z_i) \ar[r]^{\hskip 1cm \alpha_f} 
  \ar[d]_{\bigotimes\limits_{i\in I} \alpha_{g_i}}&
 \bigotimes\limits_{j\in J} \biggl(\bigotimes\limits_{k\in g^{-1}(j)} V_k(X_k)\biggr)(Z_j)
\ar[d]^{\alpha_g}
 \\
 \bigotimes\limits_{i\in I} \biggl( \bigotimes\limits_{k\in (fg)^{-1}(i)} V_k(X_k)\biggr)(Z_i) 
 \ar[r]_{\alpha_{fg}}& \bigotimes\limits_{k\in K} V_k(X_k)
 } 
 \]
 \end{prop}
 
 \noindent{\sl Proof (sketch):}  (1) Let $n=|I|$, $m=|J|$, $m_i=|f^{-1}(i)|$, so $m=\sum_{i\in I} m_i$. 
 Let also
 \[
 Z=\bigcup_{i\in I} Z_i, \quad Y=\bigcup_{j\in J} Y_j, \quad Y/i = \bigcup_{j\in f^{-1}(i)} Y_j,
 \,\, i\in I. 
 \]
 Call a $(Y,Z)$-{\em snake} a curve $S$ which is both a $Z$-snake and a $Y/i$-snake for each $i$.
 A $(Y,Z)$-snake is also a $Y$-snake. Let $\Sn(Y,Z)$ be the set of isotopy classes of $(Y,Z)$-snakes.
 It is a torsor over the subgroup $B$ in $\Br_m$ which is the wreath product of  the 
 $\Br_{m_i}$. We can view both the source and target of the desired $\alpha_f$ as pseudo-objects
 with the indexing set $\Sn(Y,Z)$, after which the map $\alpha_f$  is
 defined as
 the identity on each determination (corresponding to each $(Y,Z)$-snake).

 To prove (2), we introduce, similarly to (1), the concept of $(X,Y,Z)$-snakes 
 and  the set $\Sn(X,Y,Z)$ of their isotopy classes. 
  After this, all the arrows in the diagrams can be seen as morphisms of pseudo-objects
  indexed by $\Sn(X,Y,Z)$ and the statement becomes obvious. \qed

\paragraph{Contingency matrices.} We recall some constructions and terminology of
\cite{KS-cont}.  By a  {\em contingency matrix}  we mean a rectangular matrix 
$M=\|m_{ij}\|_{i=1,\cdots, p}^{j=1,\cdots, q}$ with $m_{ij}\in \ZZ_{\geq 0}$ 
such that each row and each column contain
at least one non-zero entry.  We will also formally include the {\empty contingency matrix}
$\emptyset$ of size $0\times 0$. 
The {\em weight} of $M$ is the number
\[
\Sigma M \,=\,\sum_{i,j} m_{ij} \,\in \, \ZZ_{>0}, \,\,\, M\neq\emptyset, \quad \Sigma\emptyset =0. 
\]
We  denote by $\CM$ the set of all contingency matrices, by
$\CM_n$ the set of all contingency matrices of weight $n$, 
by  $\CM(r,s)$ the set of all contingency matrices of size $r\times s$ and put
$\CM_n(r,s) := \CM_n\cap \CM(r,s)$. We have the {\em horizontal}  and 
{\em vertical contractions}
\[
\begin{gathered}
\del_j': \CM_n(r,s) \lra \CM_n(r,s-1), \quad i=0,\cdots, s-2,
\\
\del_i'': \CM_n(r,s) \lra  \CM_n(r-1, s), \quad j=0,\cdots. r-2,
\end{gathered}
\]
which add up the $(j+1)$st and the $(j+2)$nd column (resp. $(i+1)$st and $(i+2)$nd row). 
These contractions make the collection of $\CM_n(r,s)$ for all $r,s$ into an augmented
bi-semisimplicial  set, see \cite{KS-cont}, Prop. 1.4. 

For $M, N\in \CM_n$ we put $M\leq' N$, if $M$ can be obtained from $N$ by a series
of horizontal contractions $\del'_j$ and $M\leq '' N$, if $M$ can be obtained from $N$ by a series
of vertical contractions $\del''_i$. This defines two partial orders $\leq', \leq''$
on $\CM_n$. Thus we have the {\em elementary inequalities} $\del'_j N \leq' N$,
resp. $\del_i'' N\leq N$ and the partial order $\leq'$, resp. $\leq''$ is
generated by such elementary inequalities via transitive closure.
For any $M,L \in \CM_n$ we write
\[
\Sup(M,L) \,=\, \{ O\in\CM_n\,| M\leq' O \geq'' P\}. 
\]
An elementary inequality $\del'_j N \leq' N$ is called {\em anodyne} if, for any $j$,
among the two entries $n_{i,j+1}$ and $n_{i,j+2}$ of the $(j+1)$st and $(j+2)$nd column of $N$ that are added in $\del'_j N$, there is at least one zero. For example,
\[
\del'_0N \,=\, \begin{pmatrix} 1&2\\ 4&3\end{pmatrix} \, \leq' \, 
\begin{pmatrix} 1&0&2\\ 0&4&3\end{pmatrix} \,=\, N
\]
is anodyne. A general inequality $M\leq' N$ is called anodyne, if there is
a chain of elementary anodyne inequalities $M=M_1\leq'\cdots\leq' M_k=N$
connecting $M$ and $N$.

Similarly, an elementary inequality $\del_i'' N \leq'' N$ is called anodyne if, for each $j$.
among the two entries $n_{i+1,j}$ and $n_{i+2,j}$ there is at least one zero.
A general inequality $M\leq'' N$ is called anodyne, if there is a chain of
elementary anodyne inequalities $M=M_1\leq''\cdots\leq'' M_k=N$.

\paragraph{Components of a graded bialgebra associated to contingency matrices.}
Let $A = (A_n)_{n\geq 0}$ be a graded bialgebra in a braided monoidal
category $\Vc$. Let $M\in\CM(r,s)$ be a contingency matrix of size $r\times s$.
We put
\[
A_M \,=\,\bigotimes_{(p,q)\in \{1,\cdots, r\}\times \{1\cdots,s\}} A_{m_{p,q}}(p,q),
\]
the $2$-dimensional tensor product corresponding to $A_{m_{p,q}}$ put in the
position $(p,q)\in\RR^2$, i.e., $p+q\sqrt{-1}\in\CC$. We call $A_M$ the
{\em component} of $A$ associated to $M$ (even though it is, strictly speaking,  a tensor
product of graded components). 

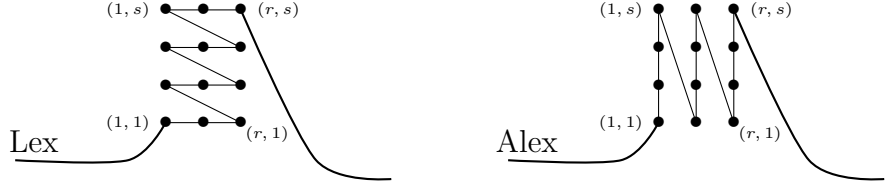
\begin{figure}[H]
\centering
\begin{tikzpicture}[scale=0.5]
\node at (1,1){\small$\bullet$}; 
\node at (1,2){\small$\bullet$}; 
\node at (1,3){\small$\bullet$}; 
\node at (1,4){\small$\bullet$}; 
\node at (2,1){\small$\bullet$}; 
\node at (2,2){\small$\bullet$}; 
\node at (2,3){\small$\bullet$}; 
\node at (2,4){\small$\bullet$}; 
\node at (3,1){\small$\bullet$}; 
\node at (3,2){\small$\bullet$}; 
\node at (3,3){\small$\bullet$}; 
\node at (3,4){\small$\bullet$}; 

\draw (1,1) -- (2,1) -- (3,1) --(1,2) -- (2,2) -- (3,2) -- (1,3) -- (2,3) -- (3,3)--
(1,4) -- (2,4) -- (3,4); 

\draw  [line width=0.8]  plot [smooth, tension=0.5]
coordinates { (-3,0) (0,0) (1,1)} ; 

\draw  [line width=0.8]  plot [smooth, tension=0.5]
coordinates { (3,4) (5,0) (7,-0.5)} ; 

\node at (-2.5, 0.5){$\Lex$}; 

\node at (0,1){\tiny$(1,1)$}; 
\node at (3.7, 0.7){\tiny$(r,1)$}; 
\node at (0,4){\tiny$(1,s)$};
\node at (4,4){\tiny$(r,s)$};

\end{tikzpicture}\quad\quad\quad
\begin{tikzpicture}[scale=0.5]
\node at (1,1){\small$\bullet$}; 
\node at (1,2){\small$\bullet$}; 
\node at (1,3){\small$\bullet$}; 
\node at (1,4){\small$\bullet$}; 
\node at (2,1){\small$\bullet$}; 
\node at (2,2){\small$\bullet$}; 
\node at (2,3){\small$\bullet$}; 
\node at (2,4){\small$\bullet$}; 
\node at (3,1){\small$\bullet$}; 
\node at (3,2){\small$\bullet$}; 
\node at (3,3){\small$\bullet$}; 
\node at (3,4){\small$\bullet$}; 

\draw (1,1) --  (1,2) -- (1,3) -- (1,4) -- (2,1) -- (2,2) -- (2,3) -- (2,4)--
(3,1) -- (3,2) -- (3,3) -- (3,4); 

\draw  [line width=0.8]  plot [smooth, tension=0.5]
coordinates { (-3,0) (0,0) (1,1)} ; 

\draw  [line width=0.8]  plot [smooth, tension=0.5]
coordinates { (3,4) (5,0) (7,-0.5)} ; 

\node at (-2.5, 0.5){$\Alex$}; 

\node at (0,1){\tiny$(1,1)$}; 
\node at (3.7, 0.7){\tiny$(r,1)$}; 
\node at (0,4){\tiny$(1,s)$};
\node at (4,4){\tiny$(r,s)$};

\end{tikzpicture}
\caption{The snakes $\Lex$ and $\Alex$ passing through
$\{1,\cdots, r\}\times\{1,\cdots, s\}$. }\label{fig:lex}
\end{figure}

 By definition,  $A_M$ is a pseudo-object with indexing set 
 $\Sn(Z)$ consisting of (isotopy classes of) snakes
 passing through $Z=\{1,\cdots, r\}\times\{1,\cdots, s\} $. 
  Among such snakes we distinguish the {\em lexicographic snake}
 $\Lex$ which reads the  elements of $Z$ one by one horizontally and the
 {\em antilexicographic snake} $\Alex$ which reads them
  one by one vertically, see  Fig. \ref{fig:lex}. The corresponding determinations of $A_M$ are
 \be\label{eq:A-M-Lex}
 \begin{gathered}
 (A_M)_\Lex \,=\, (A_{m_{1,1}}\otimes A_{m_{2,1}}\otimes\cdots \otimes
 A_{m_{r,1}}) \otimes\cdots \otimes 
 (A_{m_{1,s}}\otimes A_{m_{2,s}}\otimes\cdots \otimes
 A_{m_{r,s}}),  
 \\
  (A_M)_\Alex \,=\, (A_{m_{1,1}}\otimes A_{m_{1,2}}\otimes\cdots 
 A_{m_{1,s}}) \otimes\cdots \otimes 
 (A_{m_{r,1}}\otimes A_{m_{r,2}}\otimes\cdots \otimes
 A_{m_{r,s}}),
 \end{gathered}
 \ee
 i.e., the  ordered  tensor products of the $A_{m_{ij}}$ read along the columns, resp.
   rows of $M$. 
 
 \paragraph{Horizontal comultiplication and vertical multiplication.} 
  Our next goal is to define, for any $M\leq' N$, the {\em horizontal comultiplication map}
   $\Delta_{M,N}: A_M\to A_N$,
 and for any $M\leq'' N$, the {\em vertical multiplication map}
  $\mu_{N,M}: A_N\to A_M$,
  using the comultiplication and multiplication in $A$.

Let first
\[
\CM(r,s) \ni M=\del'_j N \,\leq' \, N\in \CM(r,s+1), \quad j\in \{1,\cdots, s-1\}
\]
 be an elementary inequality, so that for each $p$
we have:
\[
\begin{gathered}
m_{p,j+1} \,=\, n_{p, j+1} + n_{p,j+2},
\quad\quad
m_{p,q} \,=\, n_{p,q}, \quad q\leq j,
\quad\quad
m_{p,q} \,=\, n_{p, q+1}, \quad  q \geq j+2. 
\end{gathered}
\]
Introduce  closed disks $Z'_{p,q}\subset \CC $, $p=1,\cdots, r$, $q=1,\cdots, s$ by:
\[
Z'_{p,q} \,= \begin{cases}
\text{A thin ellipse encircling the points $(p, j+1)$ and $(p, j+2)$},
& \text{ if $q=j+1$};
\\
\text{the point $(p,q)$},& \text { if $q\leq j$}; 
\\
\text{the point $(p,q+1)$},& \text { if $q\geq j+2$}.  
\end{cases} 
\]
Then, on one hand, we have a canonical identification
\[
A_M \,=\,\bigotimes_{(p,q)\in \{1,\cdots, r\}\times\{1,\cdots, s\}} A_{m_{p,q}}(p,q)  \lra \bigotimes _{(p,q)\in \{1,\cdots, r\}\times\{1,\cdots, s\}}
A_{m_{p,q}}(Z'_{p,q})
\]
obtained by moving each $Z'_{p,q}$ to the point $(p,q)$  along the shortest (straight) path and then contracting it to that point if needed.
On the other hand, Proposition \ref{prop:braid-ass} gives an identification
\[
A_N\,=\,\bigotimes_{(p,q)\in \{1,\cdots, r\}\times\{1,\cdots, s\}} A'_{p,q} (Z'_{p,q}),
\]
where
\[
A'_{p,q} \,=\, \bigotimes_{(a,b)\in Z'_{p,q}} A_{n_{a,b}}\,=
 \begin{cases}
 A_{n_{p, j+1}}\otimes A_{n_{p, j+2}}, 
& \text{ if $q=j+1$};
\\
 A_{n_{p,q}},& \text { if $q\leq j$}; 
\\
 A_{n_{p, q+1}},& \text { if $q\geq j+2$}.  
\end{cases} 
\]
 Using these identifications, we define a morphism
 \[
\Delta_{M,N}: A_M\lra A_N
\]
to be the $2$-dimensional tensor 
product, over $(p,q)\in \{1,\cdots, r\}\times\{1,\cdots, s\}$,
of the morphisms $\Delta^{(p,q)}: A_{m_{p,q}}\to A'_{p,q}$ given by
\[
\Delta^{(p,q)}\,=
 \begin{cases}
 \Delta_{n_{p, j+1}, n_{p, j+2}}:  A_{m_{p, j+1}} \lra A_{n_{p, j+1}}\otimes A_{n_{p, j+2}}, 
& \text{ if $q=j+1$};
\\
\Id:  A_{m_{p,q}}\lra A_{n_{p,q}},& \text { if $q\leq j$}; 
\\
 \Id: A_{m_{p, q}}\lra A_{n_{p, q+1}},& \text { if $q\geq j+2$}
\end{cases} 
\]
and positioned at the $Z'_{p,q}$. 
In order to write $\Delta_{M,N}$   as a morphism
  of ordered tensor products without  the use of the braiding,
we can use the
  antilexicographic determinations. 
  
  \vskip .2cm
  
  Similarly, let 
  \[
 \CM(r,s) \ni M=\del''_i N \, \leq''\,  N\in \CM(r+1,s), \quad i\in\{1,\cdots, r-1\}
  \]
  be an elementary inequality, so that for
each $q$ we have
\[
m_{i+1,q}=n_{i+1,q}+ n_{i+2,q}, \quad\quad m_{p,q}=n_{p,q },\quad p\leq i,
\quad\quad m_{p,q}= n_{p+1,q},\quad p \geq i+1.
\]
 Introduce  closed disks $Z''_{p,q}\subset \CC $, $p=1,\cdots, r$, $q=1,\cdots, s$ by:
\[
Z''_{p,q} \,= \begin{cases}
\text{A thin ellipse encircling the points $(i+1,q)$ and $(i+2,q)$},
& \text{ if $p=i+1$};
\\
\text{the point $(p,q)$},& \text { if $p\leq i$}; 
\\
\text{the point $(p+1,q)$},& \text { if $p\geq i+2$}.  
\end{cases} 
\]
Then, on one hand, we have a canonical identification
\[
A_M \,=\,\bigotimes_{(p,q)\in \{1,\cdots, r\}\times\{1,\cdots, s\}} A_{m_{p,q}}(p,q)  \lra \bigotimes _{(p,q)\in \{1,\cdots, r\}\times\{1,\cdots, s\}}
A_{m_{p,q}}(Z''_{p,q})
\]
obtained by moving each $Z''_{p,q}$ to the point $(p,q)$  along the shortest (straight) path and then contracting it to that point if needed.
On the other hand,  Proposition \ref{prop:braid-ass} gives an identification
\[
A_N\,=\,\bigotimes_{(p,q)\in \{1,\cdots, r\}\times\{1,\cdots, s\}} A''_{p,q} (Z''_{p,q}),
\]
where
\[
A''_{p,q} \,=\, \bigotimes_{(a,b)\in Z''_{p,q}} A_{n_{a,b}}\,=
 \begin{cases}
 A_{n_{i+1,q}}\otimes A_{n_{i+2,q}}, 
& \text{ if $p=i+1$};
\\
 A_{n_{p,q}},& \text { if $p\leq i$}; 
\\
 A_{n_{p+1, q}},& \text { if $p\geq i+2$}.  
\end{cases} 
\]
 Using these identifications, we define a morphism
 \[
\mu_{N,M}: A_N\lra A_M
\]
to be the $2$-dimensional tensor 
product, over $(p,q)\in \{1,\cdots, r\}\times\{1,\cdots, s\}$,
of the morphisms $\mu^{(p,q)}: A''_{p,q} \to A_{m_{p,q}}$ given by
\[
\mu^{(p,q)}\,=
 \begin{cases}
 \mu_{n_{i+1, q}, n_{i+2, q}}:  A_{n_{i+1,q} }\otimes A_{n_{i+2,q}}
 \to A_{m_{i+1,q}}, 
& \text{ if $p=i+1$};
\\
\Id:  A_{m_{p,q}}\lra A_{n_{p,q}},& \text { if $p\leq i$}; 
\\
 \Id: A_{m_{p, q}}\lra A_{n_{p+1, q}},& \text { if $p\geq i+2$}
\end{cases} 
\]
and positioned at the $Z''_{p,q}$. 
In order to write $\mu_{N,M}$   as a morphism
  of ordered tensor products without  the use of the braiding,
we can use the lexicographic determinations.

\begin{prop}\label{prop:gr-bi-Hodge}
(a$'$) Let $M\leq' N$. For all chains of elementary inequalities
$M=M_1\leq'\cdots\leq' M_k=N$, 
the composition
\[
\Delta_{M,N}\,=\, \Delta_{M_{k-1}, M_k} \Delta_{M_{k-2}, M_{k-1}}\cdots \Delta_{M_1, M_2}:
A_M\lra A_N
\]
has the same value.

\vskip .2cm

(a$''$) Let $M\leq'' N$. For all chains of elementary inequalities
$M=M_1\leq'' \cdots\leq'' M_k=N$, 
the composition
\[
 \mu_{N,M}= \mu_{M_2, M_1} \mu_{M_3, M_2} \cdots \mu_{M_k, M_{k-1}}: A_N\lra A_M
\]
has the same value.

\vskip .2cm

(b) The morphisms $\Delta_{M,N}$, $M\leq' N$ and $\mu_{N,M}$, $M\leq '' N$ thus defined
satisfy the following properties:
\begin{itemize}
\item [(b1$'$)] If $L\leq' M\leq' N$, then $\Delta_{L,N}=\Delta_{M,N}\Delta_{L,M}$.

\item  [(b1$''$)] If $L\leq'' M\leq'' N$, then $\mu_{N,L}=\mu_{M,L}\mu_{N,M}$.

\item[(b2)] If $M\geq'' N \leq' L$, then
\[
\Delta_{N,L}\, \mu_{M,N} \,=\, \sum_{O\in\Sup(M,L)} \mu_{O,L} \, \Delta_{M,O}: \,\,
A_M\lra A_L. 
\]
\item[(b3$'$)] If $M\leq' N$ is an anodyne inequality, then $\Delta_{M,N}$ is an isomorphism.

\item[(b3$''$)] If $M\leq'' N$ is an anodyne inequality, then $\mu_{N,M}$ is an isomorphism. 
\end{itemize}
\end{prop}

\noindent {\sl Proof:} Parts (a$'$) and (b1$'$) follow from coassociativity of the comultiplication.
 Parts (a$''$) and (b1$''$) follow from associativity of the multiplication.
 Parts (b3$'$) and (b3$''$) are obvious from the definitions. It remains to prove (b2). 
 We do it in three steps.
 
 \vskip .2cm
 
 \noindent \ul{\sl Step 1.}  Consider first the simplest case when
 $L=(l_1, l_2)$ is a $1\times 2$ matrix, 
 $M=(m_1, m_2)^t$ is a $2\times 1$ matrix and $N=(n)$ is a $1\times 1$ matrix so that
 \[
 l_1+l_2 \,=\, m_1+m_2\,=\, n. 
 \]
In this caase $\Sup(M,L)$ consists of $2\times 2$
contingency  matrices $O=\begin{pmatrix}
o_{11}&o_{12} \\ o_{21}&o_{22}
\end{pmatrix}$ such that
\[
o_{11}+o_{21}=l_1, \,\, o_{12}+o_{22}=l_2, \quad o_{11}+o_{12}=m_1,\,\, 
o_{21}+o_{22}
=m_2. 
\]
 The claim (b2) has then the form
 \be\label{eq:b2sim}
 \Delta_{l_1, l_2} \,\mu_{m_1, m_2} \,= \sum_{O\in\Sup(M,L)}
 (\mu_{o_{11}, o_{21}}\otimes \mu_{o_{12}, o_{22}})\circ
  (\Id\otimes R_{A_{o_{12}}, A_{o_{21}}}
 \otimes\Id) \circ (\Delta_{o_{11}, o_{12}}\otimes\Delta_{o_{21}, o_{22}}), 
 \ee
 the appearance of the braiding in the middle coming from comparing the
 $\Alex$ and $\Lex$ determinations of $A_O$. 
 But this equality is simply the reformulation, at the level of graded components, 
 of the compatibility
between multiplication and comultiplication in $A$. Cf. \cite{KS-shuffle} Eq. (4.2.4). 

 \vskip .2cm
 
 \noindent \ul{\sl Step 2.} Next, suppose that  $N$ is a contingency matrix
 of arbitrary size $r\times s$ and both inequalities $M\geq'' N \leq' L$
 are elementary, so
$N=\del'_jL=\del''_i M$ for some $i$ and $j$. 
 
 \vskip .2cm
 
 The set $\Sup(M,L)$ consists then of $(r+1)\times (s+1)$
 contingency matrices $O$ such that $\del'_j O=M$, $\del''_iO=L$.
 The only part of $O$ not fixed by these conditions, is the $2\times 2$ submatrix
 $\ol O$ on rows $i+1$, $i+2$ and columns $j+1, j+2$. 
 Therefore the situation is combinatorially similar to Step 1. More precisely,
 $\Sup(M,L)$ is in bijection with the set of $2\times 2$ contingency matrices
 $\ol O= \begin{pmatrix}
 o_{i+1, j+1} & o_{i+1, j+2} \\ o_{i+2, j+1} & o_{i+2, j+2}
 \end{pmatrix}
$ such that
\[
\begin{gathered}
o_{i+1, j+1} + o_{i+1, j+2} = m_{i+1, j+1}, \quad
o_{i+2, j+1} + o_{i+2, j+2} = m_{i+2, j+1},
\\
o_{i+1, j+1} + o_{i+2, j+1} = l_{i+1, j+1}, \quad  o_{i+1, j+2} + o_{i+2, j+2}
= l_{i+1, j+2}. 
\end{gathered}
\]
Let
\[
\ol\lambda, \ol\rho: A_{m_{i+1,j+1}}\otimes A_{m_{i+2,j+1}}\lra
A_{l_{i+1, j+1}}\otimes A_{l_{i+1, j+2}}
\]
be the adaptations to our case of the  LHS and  RHS of   \eqref {eq:b2sim}, i.e.,
\[
\ol\lambda =  \Delta_{l_{i+1, j+1}, l_{i+1, j+2}}
\mu_{m_{i+1,j+1}, m_{i+2,j+1}},
\]
and $\ol\rho$ is the sum over the $\ol O$ as above. 
Thus $\ol\lambda=\ol\rho$ by Step 1. 

 \vskip .2cm

We  claim that the equality (b2) in our situation reduces to that in
Step 1, i.e., to the equality \eqref{eq:b2sim}. 
Indeed, let  $\lambda,\rho: A_M\to A_L$ be the LHS and RHS
 of  (b2). Note that each of these  morphisms  is decomposed as a $2$-dimensional tensor product
of ``elementary'' morphisms of the following $4$ types:
\begin{itemize}
\item The identity morphism from some $A_{m_{pq}}$ to some $A_{l_{p',q'}}$
for $p\neq j+1$, $q\neq i+1, i+2$, the same morphism for both $\lambda$ and $\rho$.

\item The multiplication $A_{m_{i+1,q}}\otimes A_{m_{i+2,q}}\to A_{l_{i+1,q}}$
for $q\neq i+1$, the same morphism for both $\lambda$ and $\rho$.

\item The comultiplication $A_{m_{i+1,q}}\to A_{l_{i+1,q}}\otimes A_{l_{i+2,q}}$
for $q\neq i+1, j+2$, the same morphism for both $\lambda$ and $\rho$.

\item The morphism $\ol\lambda$ for $\lambda$ and  $\ol\rho$ for $\rho$. 
\end{itemize}
This implies that $\lambda=\rho$, thus establishing Step 2.

\vskip .2cm

\begin{figure}[H] 
 \centering
   
   \begin{tikzpicture}
   
   \node at (0,0){$\bullet$};
     \node at (0,1){$\bullet$};
    \node at (1,1){$\bullet$};
      \node at (2,1){$\bullet$};
        \node at (3,1){$\bullet$};
          \node at (3,2){$\bullet$};
            \node at (4,2){$\bullet$};   
    \node at (5,2){$\bullet$};
      \node at (5,3){$\bullet$};

   \draw  (0,0) -- (5,0); 
   
   \draw (0,0) -- (0,3);

  \node at (1,0){\tiny$\bullet$}; 
    \node at (2,0){\tiny$\bullet$}; 
      \node at (3,0){\tiny$\bullet$}; 
        \node at (4,0){\tiny$\bullet$}; 
          \node at (5,0){\tiny$\bullet$};

            \node at (0,1){\tiny$\bullet$}; 
              \node at (0,2){\tiny$\bullet$}; 
                \node at (0,3){\tiny$\bullet$}; 
  
  \node at (5,3){$\bullet$};  
  \draw [dashed] (0,3) -- (5,3); 
    \draw [dashed] (0,2) -- (5,2); 
      \draw [dashed] (0,1) -- (5,1); 
 
   \draw [dashed] (1,0) -- (1,3); 
     \draw [dashed] (2,0) -- (2,3); 
       \draw [dashed] (3,0) -- (3,3); 
         \draw [dashed] (4,0) -- (4,3); 
           \draw [dashed] (5,0) -- (5,3); 
           
  \draw[ line width = 1] (0,0) -- (0,1) --(3,1) -- (3,2) -- (5,2) -- (5,3);

  \node at (-0.2, -0.2) {\tiny 0}; 
  
   \node at (1, -0.2){\tiny 1}; 
    \node at (2, -0.2){\tiny 2}; 
  \node at (3, -0.2){\tiny 3}; 
    \node at (4, -0.2){$\cdots$}; 
      \node at (5, -0.2){\tiny $a$};  

\node at (-0.2, 1) {\tiny 1}; 
\node at (-0.2, 2) {$\cdots$}; 
\node at (-0.2, 3) {\tiny $b$}; 
 
 \node at (-0.5, 0.1){$x_0$}; 
 \node at (-0.5, 1.1){$x_1$};  
\node at (5.9, 2){$x_{a+b-1}$}; 
\node at (5.7, 3){$x_{a+b}$};

   \end{tikzpicture}    
   \caption{A  taxicab path $\gamma$  in a rectangle $[0,a]\times [0,b]$.}\label{fig:taxicab}
   \end{figure}
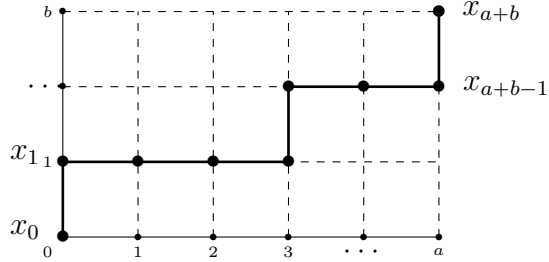

 \noindent \ul{\sl Step 3.} Let now  $M\geq'' N \leq' L$ be arbitrary. 
Let us represent
both inequalities as chains of elementary ones:
\be\label{eq:chain-M}
M=M_0 \geq'' M_1\geq'' \cdots \geq'' M_a = N \leq' M_{a+1} \leq' \cdots 
\leq' M_{a+b}=L. 
\ee
Note that   there is a unique chain  as above with given $M$ and $L$;
in particular, $N$ with $M\geq'' N \leq' L$ is unique (if it exists, which is
our assumption). 
We deduce the equality (b2)   by applying Step 2 several times.  
For this, consider {\em taxicab paths} 
\[
\gamma = [x_0, x_1, \cdots, x_{a+b}]
\]
in the rectangle $[0,a]\times [0,b]$. Such a path consists of segments
$[x_{i-1}, x_i]$, $i=1,\cdots, a+b$ of length $1$ which can be either horizontal
or vertical, with $x_0=(0,0)$ and $x_{a+b}=(a,b)$, see Fig. \ref {fig:taxicab}.

Given such $\gamma=[x_0,\cdots, x_{a+b}]$, we call a {\em $\gamma$-chain} a sequence
of contingency matrices
$L_\bullet = (L_0=M, L_1, \cdots L_{a+b}=L)$ such that:
\begin{itemize}

\item [(1)] $L_{i-1} \geq'' L_i$, if  the interval $[x_{i-1}, x_i]$ is horizontal.

\item[(2)]   $L_{i-1} \leq' L_i$, if  the interval $[x_{i-1}, x_i]$ is vertical.

\end{itemize}
Note that the equalities in a $\gamma$-chain must be elementary. 
Let $\Ch(\gamma)$ be the set of $\gamma$-chains. 
For $L_\bullet\in\Ch(\gamma)$ we have the morphism
$T_{L_\bullet}: A_M\to A_L$ defined as the composition
\[
A_M=A_{L_0} \buildrel T_{L_\bullet, 1}\over\lra A_{L_1}\buildrel T_{L_\bullet, 2}\over
\lra
\cdots \buildrel T_{L_\bullet, a+b}\over\lra A_{L_{a+b}}= A_L, 
\]
where $T_{L_\bullet, i} = \mu_{L_{i-1}, L_i}$, if the interval $[x_{i-1}, x_i]$ is horizontal
and  $T_{L_\bullet, i} = \Delta_{L_{i-1}, L_i}$, if the interval $[x_{i-1}, x_i]$ is vertical. 
The following is straightforward.

\begin{lem}
(a) For $\gamma =\gamma_{\min}:=  [(0,0), (1,0), \cdots, (a,0), (a,1), \cdots, (a,b)]$
being the  minimal (bottom right) path, there is a unique $\gamma$-chain, namely 
 \eqref{eq:chain-M}.

\vskip .2cm

(b) For $\gamma = \gamma_{\max} := [(0,0), (0,1),\cdots, (0,b), (1,b),\cdots, (a,b)]$
being the maximal (left top) path, the set $\Ch(\gamma)$ is in bijection with 
$\Sup(M,L)$. More precisely, for each $O\in Sup(M,L)$ there is a
unique $\gamma$-chain $L_\bullet$ such that $L_b=O$. \qed
\end{lem}

The lemma implies that
 \[
\begin{gathered}
\Delta_{N,L} \mu_{M,N}\,=\, \sum_{L_\bullet\in \Ch(\gamma_{\min})} T_{L_\bullet},
\quad
\sum_{O\in\Sup(M,L)}  \mu_{O,L} \, \Delta_{M,O}\,=\,
\sum_{L_\bullet\in \Ch(\gamma_{\max})} T_{L_\bullet}, 
\end{gathered}
\]
the sum in the RHS of the first equality consisting of one summand. 
So our statement reduces to the following:

\begin{lem}\label{lem:sum-inv}
The sum $\sum_{L_\bullet\in \Ch(\gamma)} T_{L_\bullet}$
is independent on the taxicab path $\gamma$ in $[0,a]\times [0,b]$. 
\end{lem} 

\noindent{\sl Proof:} It is enough to show the invariance of the sum under
an elementary modification of a path  along a $1\times 1$ square which changes 
a horizontal-then-vertical pair of unit intervals to the vertical-then-horizontal pair
completing the square. But such invariance is a consequence of Step 2,
because the inequalities corresponding to unit intervals are elementary ones. \qed

\vskip .2cm 
This  establishes Step 3 and 
Proposition \ref{prop:gr-bi-Hodge} is proved.

\vfill\eject

\section{The category of contingency matrices 
as a braided monoidal category}\label{sec:cont-mat-cat}

\paragraph{Contingency matrices as objects of a category.} We introduce a category
$\CMen$ to have, as objects, formal symbols $[M]$ for all contingency matrices
$M\in\CM$. Morphisms in $\CMen$ are generated by the generating morphisms
\[
\delta'_{M,N}: [M]\lra [N], \,\, M\leq' N, \quad \delta''_{N,M}: [N]\lra [M], \,\, M\leq'' N
\]
subject to the relations
\begin{itemize}
\item[($\CMen 1'$)] If $L\leq' M\leq' N$, then $\delta'_{L,N} = \delta'_{M,N}\delta'_{L,M}$.

\item [($\CMen 1''$)] If $L\leq'' M\leq'' N$, then $\delta''_{N,L} = \delta''_{M,L} \delta''_{N,M}$.

\item  [($\CMen 2$)] If $M\geq'' N \leq' L$, then
\[
\delta'_{N,L} \delta''_{M,N} \,=\,\sum_{O\in\Sup(M,L)} \delta''_{O,L} \delta'_{M,O}. 
\]

\item  [($\CMen 3'$)] If $M\leq' N$ is an anodyne inequality, then $\delta'_{M,N}$ is
invertible.

\item  [($\CMen 3''$)] If $M\leq'' N$ is an anodyne inequality, then $\delta''_{N,M}$ is
invertible.
\end{itemize}

More precisely, let $\CMen^+$ be the category with the objects and generating
morphisms as above which are subject to the relations ($\CMen 1'$),
($\CMen 1''$) and ($\CMen 2$). Let $E\subset \Mor_{\CMen^+}$
be the set of the $\delta'_{M,N}$, $\delta''_{N,M}$ corresponding to anodyne inequalities
$M\leq' N$, $M\leq'' N$. Then $\CMen = \CMen^+[E^{-1}]$ is the localization
of  $\CMen^+$  with respect to $\Sigma$. The set $E$ satisfies the Ore condition,
as follows from the next proposition which we leave to the reader. 

\begin{prop}\label{prop:base-change-uniq}
Let $M\geq' N \leq'' L$ be inequalities in $\CM$.

\vskip .2cm

(a) If one of these inequalities is anodyne, then $\Sup(M,L)$ consists of one element, i.e., there exists
a unique diagram of inequalities in $\CM$
 \[
\xymatrix{
O\ar[r]^{\geq''} \ar[d]_{\geq'}
& L\ar[d]^{\geq'}
\\
M\ar[r]_{\geq''}& N.
}
\]
(b) Moreover, 
  if $L\geq' N$ is anodyne, then $O\geq' M$ is anodyne.
If $M\geq'' N$ is anodyne, then $O\geq'' L$ is anodyne. \qed
\end{prop}

Let $\CMen_n\subset\CMen$ be the full subcategory on objects $[M]$, $M\in \CM_n$.
Since any inequality $M\leq' N$, $M\leq'' N$ implies equality of the weights 
$\Sigma M=\Sigma N$, the $\CMen_n$ for different $n$ are mutually orthogonal:
\[
\Hom_{\CMen}(\CMen_n, \CMen_{n'}) \,=\, 0, \quad n\neq n'. 
\]
Proposition \ref{prop:gr-bi-Hodge}  can be reformulated as follows.

\begin{cor}\label{cor:bialg-CM-V-1}
Let $A$ be a graded bialgebra in a monoidal category $\Vc$. Then the correspondence
\[
[M]\mapsto A_M, \quad \delta'_{M,N} \mapsto \Delta_{M,N}, \quad
\delta''_{N,M} \mapsto \mu_{N,M}
\]
defines a functor $\xi_A: \CMen\to\Vc$. \qed
\end{cor} 

\paragraph{Row and column exchange isomorphisms.}
Two row vectors $(m_1,\cdots, m_r)$ and $(n_1,\cdots, n_r)$ of the same size $r$
will be called {\em disjoint}, if, for each $i=1,\cdots, r$, at least one of the two numbers $m_i, n_i$
is equal to $0$, i.e.,$m_in_i=0$. Similarly for column vectors.

Let $M=\|m_{ij}\|\in\CM(r,s)$ be a contingency matrix of size $r\times s$. We denote by
\[
M_i = (m_{i1}, \cdots, m_{is}), \quad M^j=(m_{1j}, \cdots, m_{rj})^t, \quad i=1,\cdots, r, \,\,
j=1,\cdots, s,
\]
the $i$th row and the $j$th column of $M$. For $i=1,\cdots, r-1$ let $\sigma''_{i,i+1}M$
be the matrix obtained from $M$ by interchanging the $i$th and $(i+1)$st rows.
For $j=1,\cdots, s-1$ let $\sigma'_{j, j+1}M$ be the matrix obtained from $M$ by interchanging
the $j$th and $(j+1)$st columns. 

Recall that the vertical contraction $\del''_{i-1}$, $i=0, \cdots, r-1$, adds together the $i$th
and $(i+1)$st rows of a contingency matrix. 
Suppose that our $M\in \CM(r,s)$ is such that $M_i$ and $M_{i+1}$ are disjoint. Then  
the inequalities
\[
M \,\geq''\,\del''_i M = \del''_{i-1}(\sigma''_{i, i+1}M) \,\leq'' \, \sigma''_{i, i+1} M
\]
are anodyne, and we define the {\em row exchange  isomorphism}
\[
\eps''_{i, i+1}\, =\, (\delta''_{\sigma''_{i, i+1}M, \del''_{i-1} M})^{-1}\circ \delta_{M, \del''_{i-1}M}: [M] \lra [\sigma''_{i, i+1}M]
\]
in the category $\CMen$. 

Similarly, suppose that $M^j$ and $M^{j+1}$ are disjoint. Then we have anodyne inequalities
\[
M\, \geq' \del'_{j-1}M = \del'_{j-1}(\sigma'_{j, j+1}M) \,\leq' \,\sigma'_{j, j+1}M, 
\]
and we define the {\em column exchange isomorphism} in $\CMen$
\[
\eps'_{j, j+1} \, = \, \delta_{\del'_{j-1}M, \sigma'_{j, j+1}M}\circ (\delta_{\del'_{j-1}M, M})^{-1}:
[M] \lra [\sigma'_{j, j+1} M]. 
\]

\begin{prop}\label{prop:exch-braid}
(a) Let $M\in\CM(r,s)$ and $i=1, \cdots, r-2$
 be such that $M_i, M_{i+1}, M_{i+2}$ are mutually disjoint. For any permutation 
 $\tau=(\tau(1), \tau(2), \tau(3))$ of $\{1,2,3\}$ let $M_{\tau}$ be the matrix
 obtained from $M$ by permuting the $i$th, $(i+1)$st and $(i+2)$nd rows of $M$ according to
 $\tau$, e.g., $M_{(123)}=M$, $M_{(213)}=\sigma''_{i, i+1}M$ etc. Then  the  hexagon of
 row exchange isomorphisms
 \[
 \xymatrix@R=0.5em{
 & [M_{(213)}] \ar[rr]^{\eps''_{i+1, i+2}}&& [M_{(231)}] \ar[dr]^{\eps''_{i, i+1}}&
 \\
 [M_{(123)]}\ar[ur]^{\eps''_{i, i+1}} \ar[dr]_{\eps''_{i+1, i+2}}
 &&&& [M_{(321)]}
 \\
 & [M_{(132)}] \ar[rr]_{\eps''_{i, i+1}}&& [M_{(312)}] \ar[ur]_{\eps''_{i+1, i+2}}&
 }
 \]
 is commutative (braid relation).  
 
 \vskip .2cm 
 (b) A similar braid relation  for column exchange isomorphisms
  in the case when $M^j, M^{j+1}$ and $M^{j+2}$ are mutually disjoint. 
\end{prop}

\noindent{\sl Proof:} We show (a), since (b) is similar.
 By construction, each arrow in the hexagon is the composition of
two isomorphisms going through an intermediate object: one isomorphism is
 of the form $\delta''$ corresponding to an anodyne inequality $\geq''$,
the other an inverse of a $\delta''$ of this kind. Let us restore these intermediate objects
and draw the corresponding morphisms $\delta''$ (without inverting them). We get a diagram
with $12$ vertices. Let also
\[ 
\ol M \,=\,\del''_{i-1}\del''_{i-1}M \,=\,\del''_{i-1}\del''_iM 
\]
be the $(r-2)\times s$ matrix obtained by summing all three rows, the $i$th, the $(i+1)$st and
$(i+2)$nd, of $M$. Then we have an anodyne inequality  $\ol M\leq'' N$, where $N$ is any
of the $12$ matrices corresponding to the $12$ vertices of the extended diagram above. 
Therefore, putting the object $[\ol M]$ inside that diagram, we decompose it into $12$
triangles which commute because of the relation ($\CMen 1''$)  (transitivity of the maps $\delta''$).
 In this way we get
a diagram whose shape is the barycentric subdivision of the original  hexagon 
(considered as a $2$-dimensional cell complex) and
which consists of commuting triangles. This impllies the commutativity of the  ($1$-dimensional
boundary of the)
hexagon, which is the claim. \qed

\begin{prop}\label{prop:exch-diamond}
(a) Let $M\in\CM(r,s)$ and $i=1,\cdots, r-3$ have the following property:  any vector from the set
$\{M_i, M_{i+1}\}$ and any vector from the set $M_{i_2}, M_{i+3}\}$ are disjoint. For any permutation
$\ta=(\tau(1), \tau(2), \tau(3), \tau(4))$ of $\{1,2,3,4\}$ let $M_\tau$ be the matrix
obtained from $M$ by permuting the $i$th, $(i+1)$st, $(i+2)$nd and $(i+3)$rd rows according to $\tau$, e.g., $M_{(1324)}=\sigma''_{i+1, i+2}M$. Then the diagram of row exchange isomorphisms
\[
 \xymatrix@R=0.5em@C=3em{
 && [M_{(3124)}]\ar[dr]^{\eps''_{i+2, i+3}}
 &&
 \\
 [M_{(1234)}] \ar[r]^{\eps''_{i+1, i+2}} & [M_{(1324)} ] \ar[ur]^{\eps''_{i, i+1}} 
 \ar[dr]_{\eps''_{i+2, i+3}}
 &&
 [M_{(3142)}] \ar[r]^{\eps_{i+1, i+2}}& [M_{(3412)}]
 \\
 && [M_{(1342)}] \ar[ur]_{\eps''_{i, i+1}}
 }
 \]
commutes.

\vskip .2cm

(b) A similar statement for column exchange isomorphisms in the case when any vector from
$\{M^j, M^{j+1}\}$ and any vector from $\{M^{j+2}, M^{j+3}\}$ are disjoint. 

\end{prop}

\noindent {\sl Proof:} We prove (a), since (b) is similar. It suffices to prove the commutativity
of the central diamond. The argument is similar to that of Proposition \ref{prop:exch-braid}.
That is, we expand the diamond (a $4$-gon) to an $8$-gon by restoring the intermediate
objects and drawing the $\delta''$-morphisms without inverting anything. Let 
$\wt M \,=\, \del''_{i} \del''_{i-1}M_{(1324)}$ be the matrix obtained by summing the $i$th and $(i+1)$st rows
and separately summing the $(i+1)$nd and $(i+3)$rd rows of $M_{(1324)}$. Then we have an
anodyne inequality $\wt M\leq'' N$ where $N$ is any of the $8$ matrices from the $8$-gon above.
So putting $[\wt M]$ inside the $8$-gon, we fill the $8$-gon by commutative triangles which impllies that the original diamond commutes as well. \qed

\paragraph{The  monoidal structure on $\CMen$.} We make $\CMen$ into
a monoidal category by putting, on the level of objects,
\[
[M]\otimes [N] \,=\, [M\oplus N], \quad M\oplus N \,=\begin{pmatrix} M&0 \\ 0&N
\end{pmatrix}.
\]
On the level of morphisms, if $M_1\leq' M_2$ and $N_1\leq' N_2$, then
$M_1\oplus N_1\leq' M_2\oplus N_2$, and we put
\[
\delta'_{M_1, M_2}\otimes \delta'_{N_1, N_2} \,=\,\delta'_{M_1\oplus N_1, M_2\oplus N_2}.
\]
Similarly, if $M_1\leq' M_2$ and $N_1\leq' N_2$, then
$M_1\oplus N_1\leq' M_2\oplus N_2$, and we put
\[
\delta''_{M_2, M_1}\otimes\delta''_{N_2, N_1} \,=\,\delta''_{M_2\oplus N_1, M_1\oplus N_1}. 
\]
Further, let $M_1\leq' N_1$ and $M_2\leq'' N_2$. Then we have the diagram of inequalities
\[
\xymatrix{
N_1\oplus N_2 \ar[r]^{\geq''} 
\ar[d]_{\geq'} &
 N_1\oplus M_2
 \ar[d]^{\geq '}
\\
M_1\oplus N_2 \ar[r]_{\geq''} & M_1\oplus M_2
}
\]
and $\Sup(M_1\oplus N_2, N_1\oplus M_2) = \{N_1\oplus N_2\}$ consists of one element.
Therefore
\[
\delta'_{M_1\oplus M_2, N_1\oplus M_2} \delta''_{M_1\oplus N_2, M_1\oplus M_2} \,=\,
\delta''_{N_1\oplus N_2, N_1\oplus M_2} \delta' _{M_1\oplus N_2, N_1\oplus N_2}
\]
and we define $\delta'_{M_1, N_1}\otimes\delta''_{N_2, M_2}$ to be equal to this common value.

\vskip .2cm

Similarly, let $M_1\leq'' N_1$ and $M_2\leq' N_2$. We have the diagram of inequalities
 \[
\xymatrix{
N_1\oplus N_2 \ar[r]^{\geq''} 
\ar[d]_{\geq'} &
 M_1\oplus N_2
 \ar[d]^{\geq '}
\\
N_1\oplus M_2 \ar[r]_{\geq''} & M_1\oplus M_2
}
\]
and put
\[
\delta''_{N_1, M_1}\otimes \delta'_{M_2, N_2} \,:=\,
\delta'_{M_1\oplus M_2, M_1\oplus N_2} \delta'' _{N_1\oplus M_2, M_1\oplus M_2} \,=\,
\delta''_{N_1\oplus N_2, M_1\oplus N_2} \delta'_{N_1\oplus M_2, N_1\oplus N_2},
\]
the second inequality following from
$\Sup(N_1\oplus M_2, M_1\oplus N_2) = \{N_1\oplus N_2\}$. 

\begin{prop}
The above data on objects and generating morphisms define a monoidal structure $\otimes$
on $\CMen$ with unit object $\1=[\emptyset]$. 
\end{prop}

\noindent{\sl Proof:} By construction, the operation $\otimes$ is strictly associative on objects. 
What remains
to prove is that $\otimes$ extends to a functor in each argument, i.e., that our definitions
are compatible with the relations in $\CMen$. 
For this, we proceed as follows.

First, our definitions imply that for any two generating morphisms
$f: [M_1]\to [M_2]$, $g: [N_1]\to [N_2]$ we have
\[
f\otimes g\,=\, (f\otimes \Id_{[N_2]})(\Id_{[M_1]}\otimes g) \,=\,
(\Id_{[M_2]}\otimes g)(f\otimes\Id_{[N_1]}). 
\]
So it suffices to show that for any $M,N\in\CM$  the operations
$(-\otimes\Id_{[N]})$, $(\Id_{[M]}\otimes -)$ on generating morphisms
preserve the relations in $\CMen$. We consider $(-\otimes \Id_{[N]})$,
the case of $(\Id_{[M]}\otimes -)$ being similar. 

\vskip .2cm

For ($\CMen 1'$),  ($\CMen 1'$) such preservation is obvious.
For  ($\CMen 2$) it follows from the identification
\[
\Sup(L,M) \lra \Sup(L\+ N, M\+ N), \quad O\mapsto O\+ N. 
\]
For  ($\CMen 3'$),  ($\CMen 3''$) it follows from the following obvious fact:
if $M_1\leq' M_2$ is anodyne, then $M_1\+ N\leq' M_2\+ N$ is anodyne
also, and similarly for $M_1\leq'' M_2$. \qed

\paragraph{Braiding on $\CMen$.} Let $M\in\CM(p,q)$ and $N\in\CM(r,s)$.
We define the braiding isomorphism
\[
R_{[M], [N]}:  [M] \x [N] \,=\, [M\+ N] \lra [N\+ M] \,=\, [N]\x [M]
\]
by mimicking the standard Eckmann-Hilton procedure  in topology
  (``jeu de taquin'' proving the  commutativity of $\pi_2$).
 More precisely, we define $R_{[M], [N]}$ as the composition
 \[
 [M\oplus N] =
 \left[\begin{pmatrix} M&0\\ 0&N
\end{pmatrix}\right]
 \buildrel R''_{M,N}\over\lra
 \left[  \begin{pmatrix} 0&N\\ M&0
\end{pmatrix}\right] \buildrel R'_{M,N}\over\lra 
\left[
 \begin{pmatrix} N&0\\ 0&M
\end{pmatrix}\right]
=[N\oplus M], 
 \]
 where: 
 \begin{itemize}
 \item $R''_{M,N}$ is the composition of $pr$ row exchange isomorphisms
 moving $r$ rows of $(0 N)$ past the $p$ rows of $(M 0)$. This can be done
 in several ways but Proposition
 \ref{prop:exch-diamond}(a) mplies that all of them lead to the same result, which is denoted $R''_{M,N}$.

 \item $R'_{M,N}$ is the composition of $qs$ column exchange isomorphisms moving $s$ columns of
 $\begin{pmatrix} N\\0\end{pmatrix}$ past $q$ columns of 
  $\begin{pmatrix} 0\\M\end{pmatrix}$. Again, this can be done in several ways but
  Proposition
 \ref{prop:exch-diamond}(b) mplies that all of them lead to the same result, which is denoted $R'_{M,N}$. 
 \end{itemize}

\begin{prop}\label{prop:CM-braided}
The isomorphisms $R_{[M], [N]}$ make $\CMen$ into a braided monoidal
category. 
\end{prop}

\noindent{\sl Proof:} We first show that the  $R_{[M], [N]}$  are natural in each variable. Naturality in the first variable
means that for any morphism $\phi: [M]\to [L]$  and any object $[N]$ in $\CMen$ the diagram
(the naturality square)
\[
\xymatrix{
[M]\otimes [N] \ar[d]_{\phi\otimes [N]}
\ar[rr]^{R_{[M], [N]}} && [N]\otimes [M]  \ar[d]^{[N[\otimes \phi} 
\\
 [L]\otimes [N] \ar[rr]_{R_{[L], [N]}}&& [N]\otimes [L]
}
\]
is commutative.   
To show this, it suffices to assume that $\phi$ is one of the elementary generating  morphisms,
i.e., we are in either of the two cases:
\begin{itemize}
\item[(i)] $\phi=\delta'_{ML}$, where $M=\del'_jL$ is obtained from $L$ by  a horizontal contraction (adding two adjacent columns);

\item[(ii)] $\phi=\delta''_{ML}$, where $L=\del''_iM$ is obtained from $M$ by a 
vertical contraction (adding two adjacent rows). 
\end{itemize}
Consider the case (i). The naturality square whose commutativity we need to prove,
decomposes into two:
\be\label{eq:two-sq}
\xymatrix{
 \left[ {\begin{pmatrix} M & 0 \\ 0&N \end{pmatrix} }\right] \ar[r]^{R''_{MN}} 
 \ar[d]_{\phi\otimes[N]}
 & 
  \left[ {\begin{pmatrix} 0&N \\ M&0 \end{pmatrix} }\right] \ar[r]^{R'_{MN}}
  \ar[d]^\psi
  & 
   \left[ {\begin{pmatrix} N & 0 \\ 0&M \end{pmatrix} }\right]\ar[d]^{[N]\otimes\phi}
   \\
 \left[ {\begin{pmatrix} L & 0 \\ 0&N \end{pmatrix} }\right] \ar[r]^{R''_{LN}} & 
  \left[ {\begin{pmatrix} 0&N \\ L&0 \end{pmatrix} }\right] \ar[r]^{R'_{LN}}& 
   \left[ {\begin{pmatrix} N & 0 \\ 0&L \end{pmatrix} }\right],  
}
\ee
 where $\psi$ is induced by the inequality $ \begin{pmatrix} 0&N \\ M&0 \end{pmatrix} \leq' 
  \begin{pmatrix} 0&N \\ L&0 \end{pmatrix}$. Note that the other two vertical arrows are,
  by construction, also induced by the corresponding inequalities. We prove the commutativity
  of each of the two squares separately. 
  
  \vskip .2cm
  
  \noindent \ul{Left square:} The morphisms $R''_{MN}$ and $R''_{LN}$ are compositions of row
  exchange isomorphisms, i.e., of  $\delta''$-isomorphisms induced by anodyne vertical contractions
  and of the inverses of such isomorphisms. These isomorphisms go through intermediate
  objects corresponding to  matrices  obtained from ${\begin{pmatrix} M & 0 \\ 0&N \end{pmatrix} }$
  and ${\begin{pmatrix} L & 0 \\ 0&N \end{pmatrix} }$ by some number  of row exchanges and then, possibly,
   summation of two disjoint adjacent rows. Let us restore these intermediate objects
   and the $\delta''$-isomorphisms connecting them, without inverting these isomorphisms. In this way
   we replace  the square  by a  diagram  of the form
    \be\label{eq:ladder1}
   \xymatrix{
    \left[{\begin{pmatrix} M & 0 \\ 0&N \end{pmatrix}}\right] 
    \ar[d]_{\phi\otimes[N]}^{\delta'}
    \ar[r]^{\hskip 1cm \delta''} &\bullet  \ar@{-->} [d]_{\delta'}
     &\ar[l]_{\delta''}\bullet  \ar@{-->} [d]_{\delta'}
      \ar[r]^{\delta''}&\bullet  \ar@{-->} [d]_{\delta'}
    &\ar[l]_{\delta''} \cdots \ar[r]^{\delta''}&\bullet  \ar@{-->} [d]_{\delta'}
     &\ar[l]_{\hskip -1cm \delta''}
      \left[ {\begin{pmatrix} 0&N \\ M&0 \end{pmatrix} }\right] 
       \ar[d]^\psi_{\delta'}
\\
 \left[{\begin{pmatrix} L & 0 \\ 0&N \end{pmatrix}}\right] 
  \ar[r]_{\hskip 1cm \delta''} &\bullet &\ar[l]^{\delta''}\bullet \ar[r]_{\delta''}&\bullet
    &\ar[l]^{\delta''} \cdots \ar[r]_{\delta''}&\bullet &\ar[l]^{\hskip -1cm \delta''}
      \left[ {\begin{pmatrix} 0&N \\ L&0 \end{pmatrix} }\right] 
     }
   \ee
    where the horizontal maps are $\delta''$-isomorphisms. Note further that we have morphisms 
  between the corresponding intermediate objects indicated by the dotted vertical arrows. They correspond to
  the horizontal contractions of the intermediate matrices. We obtain a ladder diagram consisting of many squares,
  with horizontal maps being $\delta''$-isomorphisms and vertical maps being $\delta'$-morphsms. 
  We claim that each of these squares is commutative. Indeed, such a square corresponds to a  square of inequalities of
  the type discussed in Proposition \ref{prop:base-change-uniq}: two of the inequalities  of the same type 
  (in our case, $\leq''$)
  are anodyne, 
  In this situation Proposition \ref{prop:base-change-uniq} and the relation  ($\CMen 2$)  give that the square is commutative,
  as the sum in  ($\CMen 2$) consists of one summand. This implies that the boundary  of the entire diagram, 
   formed by 
  inverting the  isomorphisms oriented $\leftarrow$, i.e., the left square in \eqref{eq:two-sq}, is commutative. 
  
 \vskip .2cm
 
 \noindent \ul{Right square:} The morphisms $R'_{MN}$ and $R'_{LN}$ are compositions of column
  exchange isomorphisms, i.e., of  $\delta'$-isomorphisms induced by anodyne horizontal contractions
  and of the inverses of such isomorphisms. Restoring the itnermedaite objects involved in these
  isomorphisms, we obtain a diagram somewhat similar to \eqref{eq:ladder1}:
    \[
   \xymatrix{
   \left[{\begin{pmatrix} 0&N \\ M&0  \end{pmatrix}}\right] 
    \ar[d]_{\psi}^{\delta'}
    &\ar[l]_{\hskip 1cm \delta'} \bullet  \ar[r]^{\delta'} & \bullet   & \ar[l] 
    _{\delta''} \bullet\ar[r] ^{\delta'}   & \cdots
   & \ar[l] _{\delta'} \bullet\ar[r]^{ \hskip -1cm \delta'}   & 
     \left[ {\begin{pmatrix} N&0 \\ 0&M \end{pmatrix} }\right] \ar[d]^{[N]\otimes \phi}_{\delta'}
     \\
       \left[{\begin{pmatrix} 0&N \\ L&0 \end{pmatrix}}\right]  &\ar[l]^{\hskip 1cm \delta'}  \bullet \ar[r]_{\delta'} & \bullet& \ar[l]^{\delta'}  \bullet\ar[r]_{\delta'} & \cdots
   & \ar[l]^{\delta'}  \bullet\ar[r]_{\hskip -1cm \delta'}& 
     \left[ {\begin{pmatrix} N&0 \\ 0&L \end{pmatrix} }\right]. 
     }
   \]
 This diagram consists entirely of $\delta'$-morphisms. Further, unlike  \eqref{eq:ladder1},
 the bottom row here is longer than the top one, 
 since $M=\del'_jL$ has one fewer column than $L$,
 being obtained from $L$ by adding the $(j+1)$st and $(j+2)$nd columns. Let us denote
 these columns for short by $\bl=L^{j+1}$ and $\bl'=L^{j+2}$. 
 
 Now, some objects in the bottom row can be assigned  ``matches'' in the top one, from which they  receive   $\delta'$-maps which we  add to the diagram  as vertical arrows.
  These objects correspond to matrices
  which contain the columns $\begin{pmatrix} 0\\  \bl \end{pmatrix}$ and 
 $\begin{pmatrix} 0\\  \bl' \end{pmatrix}$  
 situated next to each other, and the corresponding matching
 matrix in the top row is obtained by adding these columns. 
 In this way we get several vertical arrows which decompose our diagram into fragments
 of two types.
 
 A fragment of the first type is a square obtained when two vertical arrows are positioned next to
 each other. Each such square is commutative by transitivity of $\delta'$-maps.
 
 A fragment of the second type is obtained when a column of $N$, denote it $\bn$,
 is moved past $\bl$ and $\bl'$. Such a fragment has the form
 (we do not depict any other columns that are unchanged throughout the procedure):
  \[
 \xymatrix{
  \left[ {\begin{pmatrix} 0&\bn \\ \bl + \bl' & 0 \end{pmatrix}}\right]
  \ar[d]_{\delta'} && \ar[ll]_{\delta'} 
   \left[ {\begin{pmatrix} \bn \\ \bl + \bl'  \end{pmatrix}}\right]
      \ar@{-->} [dl]_{\delta'}  \ar@{-->} [d]_{\delta'}  \ar@{-->} [dr]^{\delta'}
\ar[rr]^{\delta'} &&
    \left[ {\begin{pmatrix} \bn & 0 \\ 0 & \bl + \bl'  \end{pmatrix}}\right]
    \ar[d]^{\delta'}
\\
 \left[ {\begin{pmatrix}0&0&\bn \\ \bl&\bl'&0\end{pmatrix}}\right]& \ar[l]^{\hskip 0.5cm \delta'} 
  \left[ {\begin{pmatrix}0&\bn \\ \bl & \bl' \end{pmatrix}}\right] \ar[r]_{\hskip -0.5cm \delta'} &
    \left[ {\begin{pmatrix} 0&\bn& 0 \\ \bl & 0 & \bl'  \end{pmatrix}}\right] &
    \ar[l]^{\hskip 0.5cm \delta'}   \left[ {\begin{pmatrix} \bn & 0 \\ \bl & \bl'  \end{pmatrix}}\right]
    \ar[r]_{\hskip -0.5cm \delta'} & 
      \left[ {\begin{pmatrix} \bn & 0 & 0 \\ 0 & \bl & \bl'   \end{pmatrix}}\right]
 }
 \]
To show that this fragment becomes commutative after inverting 
the arrows oriented $\leftarrow$,
we decompose it by the dotted arrows (which are
likewise $\delta'$-morphisms) into two $4$-gons
and two triangles. Each of them is commutative by transitivity of $\delta'$-morphisms. 

This proves that the right square in \eqref{eq:two-sq} is commutative in the situation of
 Case (i) above, i.e,. under the assumption that $\phi=\delta'_{ML}$, where $M=\del'_jL$.
 In this way we show the naturality of the $R_{[M], [N]}$ in  the first variable in Case (i).
 
 Naturality in the first argument in
Case (ii) when   $\phi=\delta''_{ML}$, $L=\del''_iM$, is analyzed completely analogously
 except the roles of the left and right squares in  \eqref{eq:two-sq}  will be interchanged. 
 
 Further, the naturality in the second argument is  also completely analogous. 
 This proves that the $R_{[M], [N]}$ are natural in both arguments. 
 
 \vskip .2cm
 
 To prove that $R$ is a braiding, it remains to show the commutativity of the {\em braiding
 triangles} \cite{bakalov, joyal-street}. These triangles are of two classes. The triangles of the first class have the form
 \[
 \xymatrix{
 [M]\otimes [M'] \otimes [N] \ar[d]_{R_{[M]\otimes [M'], [N]}}
 \ar[rr]^{[M] \otimes R_{[M'],[N]}}&& [M]\otimes [N] \otimes [M']
 \ar[dll]^{ \hskip .5cm R_{[M], [N]}\otimes [M']}
 \\
 [N]\otimes [M] \otimes [M']
 }
 \]
 for any three objects $[M], [M'], [N]$. The triangles of the second class are similarly
 associated to any $[M], [N], [N']$ and express two ways of passing from
 $[M]\otimes [N]\otimes [N']$ to $[N]\otimes [N']\otimes [M]$. 
 The commutativity of such triangles follows straightforwardly from Propositions
 \ref{prop:exch-braid} (braid relation for row or column exhanges) and 
 \ref{prop:exch-diamond}. 
 Proposition \ref{prop:CM-braided} is proved. 

\vskip .2cm

 We now notice the following refinement of Corollary \ref{cor:bialg-CM-V-1}. 
 
 \begin{prop}\label{prop:bialg-CM-V}
 In the situation of  Corollary \ref{cor:bialg-CM-V-1}, the functor $\xi_A: \CMen\to \Vc$
 is braided monoidal. 
 \end{prop}
 
 \noindent{\sl Proof:} We first construct isomorphisms
 \[
 \xi_A([M])\otimes\xi_A([N]) = A_M\otimes A_N \buildrel \phi_{M,N} \over\lra
 A_{M\oplus N} = \xi_A([M]\otimes [N]).
 \]
 Suppose $M$ is of size $p\times q$ and $N$ is of size $r\times s$. By definition, the component
  $A_M$, being a $2$-dimensional tensor product,  is a  pseudo-objecr and as such, is given in terms of determinations corresponding to snakes. 
 In particular
 \eqref {eq:A-M-Lex}, the determination  $(A_M)_{\Lex}$ corresponding to the $\Lex$ snake, 
 is  the ordered tensor
 product of the $A_{m_{ij}}$ along the columns of $M$. Similarly, $(A_M)_{\Alex}$, the determination
 corresponding to the $\Alex$ snake, is the ordered tensor
 product of the $A_{m_{ij}}$ along the rows of $M$. They are idenfitied by the braiding $R_{b_M}$,
 where $b_M\in \Br_{pq}$ is the braid (depending only on $p$ and $q$) connecting $\Lex$ and $\Alex$ snakes for $M$ (in fact, for any $p\times q$ matrix).
   Similarly for $(A_N)_\Lex$ and $(A_N)_\Alex$ and $R_b$ for  $b_N\in\Br_{rs}$. 
 
 \vskip .2cm
 
 Now note that reading $M\oplus N$ along the columns (and ignoring the $0$'s in the off-diagonal blocks) is the same as first reading $M$ along the columns and then reading $N$ in the same way.
  This gives an isomorphism
 $\phi_{M,N}^\Lex: (A_M)_\Lex\otimes (A_N)_\Lex\to (A_{M\oplus N})_\Lex$. 
 Similarly, reading $M\oplus N$ along the rows (and ignoring the $0$s as above) is the same as
 first reading $M$ and then reading $N$ in this way. This gives an isomorphism
 $\phi_{M,N}^\Alex: (A_M)_\Lex\otimes (A_N)_\Alex\to (A_{M\oplus N})_\Alex$. We claim that 
  $\phi_{M,N}^\Lex$ and  $\phi_{M,N}^\Alex$ give the same morphism of pseudo-objects
  $\phi_{M,N}: A_M\otimes A_N\to A_{M\oplus N}$. Indeed, consider the juxtaposition (direct sum)
  homomorphism
  \[
  \oplus: \Br_{pq} \times \Br_{rs} \to \Br_{pq+rs}. 
  \]
  The $\Lex$ and $\Alex$ determinations of $A_{M\oplus N}$ are related by the braiding
  $R_{b_{M\oplus N}}$, where $b_{M\oplus N}\in\Br_{pq+rs}$ is the braid relating the $\Lex$ and $\Alex$
  snakes for block-diagonal matrices $(p+r)\times (q+s)$. We notice that
  $b_{M\oplus N} = b_M\oplus b_N$ and therefore we have a commutative square
  \[
  \xymatrix{
  (A_M)_\Lex\otimes (A_N)_\Lex
  \ar[d]_{ R_{b_M}\otimes R_{b_N}}
   \ar[r]^{\hskip 0.5cm \phi_{M,N}^\Lex} & (A_{M\oplus N})_\Lex
   \ar[d]^{R_{b_{M\oplus N}}}
   \\
    (A_M)_\Alex\otimes (A_N)_\Alex \ar[r]^{\hskip 0.5cm \phi_{M,N}^\Lex} & (A_{M\oplus N})_\Alex
  }
  \]
  which implies that the resulting morphism of pseudo-objects is the same for both $\Lex$ and
  $\Alex$ determinations. This defines
  $\phi_{M,N}$. 
  
  \vskip .2cm
  
  Next, we show that the $\phi_{M,N}$ are natural in $[M]$ and $[N]$. It suffice to check the naturality on
  generating morphisms $\delta'$ or $\delta''$ for $M$ or $N$. Naturality for $\delta'$ (horizontal contractions,
  adding some adjacent columns) is immediate in the $\Alex$ determination. Indeed, when reading the $A_{m_{ij}}$
  along the rows,  the action of $\delta'$, i.e., horizontal comultiplication, will respect the order of the product,
  i.e., will produce new tensor factors in positions which are adjacent with respect to the order. 
  Similarly, naturality for $\delta''$ is immediate in the $\Lex$ determination, reading the
  $A_{m_{ij}}$ along the columns. 
  
  This naturality makes $\xi_A$ into a monoidal functor. It remains to show that $\xi_A$
  is in fact a braided monoidal functor, i.e., preserves the braiding. This  verification 
  is straightforward and left
  to the reader. Proposition \ref{prop:bialg-CM-V} is proved.

\vfill\eject

\section{The category of contingency matrices and the PROB $\Ben$ of graded
bialgebras}\label{sect:CM-Ben}

\paragraph{The graded bialgebra $\aen$ in $\CMen$.} We now define
a graded bialgebra $\aen=(\aen_n)_{n\geq 0}$ in $\CMen$
with components $\aen_n=[n]$ (the object corresponding to the $1\times 1$
contingency matrix $(n)$) for $n>0$ and $\aen_0=\1=[\emptyset]$. 
The multiplication and comultiplication are given by
\[
\begin{gathered}
\mu_{m,n}: [m]\otimes [n] \,=\left[\begin{pmatrix} m&0\\0&n
\end{pmatrix} \right] \buildrel (\delta')^{-1}\over\lra  \left[
\begin{pmatrix} m\\ n\end {pmatrix}\right] \buildrel \delta''\over\lra
[m+n],
\\
\Delta_{m,n}: [m+n] \buildrel \delta'\over\lra [(m, \, n)] 
\buildrel (\delta'')^{-1}\over\lra \left[\begin{pmatrix} m&0\\0&n
\end{pmatrix} \right] \,=\, [m]\otimes [n]. 
\end{gathered}
\]

\begin{prop}\label{prop:aen-in-CM}
The morphisms $\mu_{m,n}$, $\Delta_{m,n}$ make $\aen$ into a graded
bialgebra in $\CMen$. 
\end{prop}

\noindent{\sl Proof:} We first prove associativity. For this, we must compare two morphisms
$[m]\otimes [n]\otimes [p]\to [m+n+p]$ corresponding to two bracketing of the triple product. There
morphisms are the compositions of the upper and lower paths in the boundary of the 
following diagram, the paths obtained by inverting the $\delta'$-isomorphisms:
\[
\xymatrix{
&
\ar[dl]_{\delta'}
 \left[ {\begin{pmatrix}
 m&0 \\n&0\\0&p 
 \end{pmatrix}}\right]
\ar[r]^{\delta''} & 
 \left[ {\begin{pmatrix}
 m+n&0 \\ 0&p
  \end{pmatrix}}\right]
  & \ar[l]_{\hskip .5cm \delta'} 
   \left[ {\begin{pmatrix} m+n \\ p
    \end{pmatrix}}\right]
    \ar[dr]^{\delta''}
    &
    \\
   \left[ {\begin{pmatrix}
   m&0&0 \\ 0&n&0 \\ 0&0&p
   \end{pmatrix}}\right]   
   &
    \left[ {\begin{pmatrix}
    m\\ n\\ p
      \end{pmatrix}}\right]  
      \ar@{-->}[u]^{\delta'}  \ar@{-->}[d]_{\delta'}
       \ar@{-->}[urr]^{\delta''}  \ar@{-->}[drr]_{\delta''}
   &&&
   [m+n+p]
   \\
   & 
   \ar[ul]^{\delta'}
      \left[ {\begin{pmatrix}
      m&0\\0&n\\0&p
     \end{pmatrix}}\right]   
     \ar[r]_{\delta''} & 
    \left[ {\begin{pmatrix}     
 m&0\\0&n+p  
     \end{pmatrix}}\right]   
     & 
     \ar[l]^{\delta'} 
     \left[ {\begin{pmatrix}
m\\ n+p   
  \end{pmatrix}}\right]. 
  \ar[ur]_{\delta''} 
}
\]
To prove that these two paths have the same composition, we decompose the diagram into four $4$-gons
by the dotted arrows as shown and notice that each of these $4$-gons is commutative.

Indeed, the leftmost $4$-gon commutes by transitivity of $\delta'$-morphisms. The rightmost 
$4$-gon commutes by transitivity of $\delta''$-morphisms. The remaining two $4$-gons commute by  the relation
 ($\CMen 2$) since the sum in that  relation consists of one summand by Proposition 
 \ref{prop:base-change-uniq}. 
 
 This proves associativity of $\mu$. The proof of coassociativity  of $\Delta$ is similar. 
 
 Finally, we prove compatibility of $\Delta$ and $\mu$. This is expressed by Eq. \eqref{eq:b2sim},
 we we assume that we are in the situation of \eqref{eq:b2sim}. The composition
 $\Delta_{l_1, l_2} \,\mu_{m_1, m_2} $ is, in our case, given by the border (top horizontal followed
 by the right vertical) path in the following diagram:
 \[
 \xymatrix{
  \left[ {\begin{pmatrix} 
  m_1&0\\ 0&m_2
    \end{pmatrix}}\right]
    \ar[r]^{ \hskip 0.5cm (\delta')^{-1} }_{\hskip 0.5cm \beta} &
   \left[ {\begin{pmatrix} 
   m_1 \\ m_2
    \end{pmatrix}}\right]
    \ar[d]_{\delta'}^{\psi_O}
     \ar[r]^{ \hskip 0.5cm \delta''}_{ \hskip 0.5cm \phi }& 
    [n]   \ar[d]^{\delta'}_\psi 
    \\
    &  
      \left[ {\begin{pmatrix} 
      o_{11}&o_{12} \\ o_{21} & o_{22}
    \end{pmatrix}}\right] 
    \ar[r]_{ \hskip 0.5cm \phi_O}^{ \hskip 0.5cm \delta''} & [(l_1, l_2)]  
    \ar[d]^{(\delta'')^{-1}}_{\alpha}   
    \\
    &&
     \left[ {\begin{pmatrix} 
     l_1&0 \\ 0&l_2
    \end{pmatrix}}\right]       
 }
 \]
 Here the matrix $O= {\begin{pmatrix} 
      o_{11}&o_{12} \\ o_{21} & o_{22}
    \end{pmatrix}}$ is a (so far arbitrary)  element of $\Sup(M,L)$. We denoted  for short by
    $\psi, \psi_O, \phi, \phi_O$ the $\delta'$ and $\delta''$-morphisms in the square and by $\alpha$and $\beta$
    the inverted isomorphisms at the end and the beginning of the border path. 
    By the relation ($\CMen 2$) we have 
    \[
    \psi \phi \,=\sum_{O\in\Sup(M,L)} \phi_O \psi_O. 
    \]
 So it suffices to show that for each $O\in\Sup(M,L)$ we have
 \be\label{eq:alpha-phi-psi-beta}
 \alpha \phi_O \psi_O\beta \,=\, (\mu_{o_{11}, o_{21}}\otimes \mu_{o_{12}, o_{22}})\circ
  (\Id\otimes R_{A_{o_{12}}, A_{o_{21}}}
 \otimes\Id) \circ (\Delta_{o_{11}, o_{12}}\otimes\Delta_{o_{21}, o_{22}}). 
 \ee
 so that the summands in  the RHS of \eqref{eq:b2sim}  match those in ($\CMen2$). 
 We represent the two sides of \eqref{eq:alpha-phi-psi-beta} by the upper and lower path
 in the
 boundary of the following diagram (more precisely, the paths, going from left to right, are obtained
 by inverting the isomorphisms oriented the other way):
 {\tiny
 
 \[
 \xymatrix@C=0.7em{
 &
   \left[ {\begin{pmatrix} 
    m_1\\m_2
       \end{pmatrix}}\right] 
       \ar[dl]_{\delta'}^{\sim}
         \ar[r]^{\hskip -0.5cm \delta'}_{\hskip -0.5cm \psi_O} & 
          \left[ {\begin{pmatrix} 
      o_{11}&o_{12} \\ o_{21} & o_{22}
    \end{pmatrix}}\right] 
     \ar@{-->}[ddll]_{\kappa}   \ar@{-->}[ddrr]^\rho
    \ar[r]^{  \hskip 0.3cm \delta''}_{ \hskip 0.3cm \phi_O} & [(l_1, l_2)] &
    \\
  \left[ {\begin{pmatrix} 
  m_1&0\\ 0&m_2
    \end{pmatrix}}\right]
    \ar[d]_{\delta'}
     &&&&   
         \left[ {\begin{pmatrix} 
     l_1&0 \\ 0&l_2
    \end{pmatrix}}\right]   
    \ar[d]^{ \delta''}
    \ar[ul]_{\delta''}^\sim 
    \\
   \left[ {\begin{pmatrix} 
   o_{11} & o_{12} &0&0 \\ 0&0&o_{21}&o_{22}
    \end{pmatrix}}\right]     
    &&&&  
       \left[ {\begin{pmatrix} 
       o_{11}&0 \\ o_{21} & 0 \\ 0&o_{12} \\ 0&o_{22}   
           \end{pmatrix}}\right] 
           \ar[dl]^{\delta'}  
           \\
  &    
    \left[ {\begin{pmatrix} 
    o_{11}&0&0&0 \\
    0&o_{12}&0&0 \\
    0&0& o_{21}&0\\
    0&0&0&o_{22}    
     \end{pmatrix}}\right]  
     \ar[ul]^{\delta''}
     \ar[rr]_{\Id\otimes R\otimes Id} && 
       \left[ {\begin{pmatrix} 
    o_{11}&0&0&0 \\
    0&o_{12}&0&0 \\
    0&0& o_{21}&0\\
    0&0&0&o_{22}    
     \end{pmatrix}}\right]  
     &
        }
 \]
 }  
 To show the equality of the compositions of these paths, we decompose the diagram into two $4$-gons
 and a pentagon
  by the dotted arrows
 $\kappa$ and $\rho$, where:
 \begin{itemize}
 \item $\kappa$ is the composition
 \[
         \left[ {\begin{pmatrix} 
      o_{11}&o_{12} \\ o_{21} & o_{22}
    \end{pmatrix}}\right]  \buildrel\delta'\over\lra 
    \left[ {\begin{pmatrix} 
    o_{11}&0&o_{12}&0 \\  
    0&o_{21}&0&o_{22}
       \end{pmatrix}}\right]  
       \buildrel {\text{Column exchange}}\over\lra 
        \left[ {\begin{pmatrix} 
        o_{11}& o_{12}&0&0 \\
        0&0&o_{21}&o_{22}
              \end{pmatrix}}\right],  
 \]
 so it composed entirely of $\delta'$-morphisms and their inverses. 
 
 \item $\rho$  is the composition
  \[
         \left[ {\begin{pmatrix} 
      o_{11}&o_{12} \\ o_{21} & o_{22}
    \end{pmatrix}}\right] 
     \buildrel (\delta'')^{-1}\over\lra 
      \left[ {\begin{pmatrix} 
      o_{11}&0\\
      o&o_{12}\\
      o_{21}&0\\
      0&o_{22}
              \end{pmatrix}}\right] 
            \buildrel {\text{Row exchange}}\over\lra 
               \left[ {\begin{pmatrix} 
o_{11}&0\\
o_{21}&0\\
0&o_{12}\\
0&o_{22}
   \end{pmatrix}}\right], 
   \]
 \end{itemize}
 so it is composed entirely of $\delta''$-morphisms and their inverses. 
 The left $4$-gon in the decomposed diagram is commutative by transitivity of $\delta'$-morphisms.
 The right $4$-gon is commutative by transitivity of $\delta''$-morphsms.
 Finally, the pentagon consists  entirely of  anodyne $\delta'$- or $\delta''$-isomorphisms 
 and their inverses which  move the $o_{ij}$ around in the plane. We can view them as moving $4$ points
 in the plane. After we go around the pentagon, we return to the same position. Moreover,
 the braid on $4$ strands representing this move, is trivial. This triviality of the braid
 implies the commutativity of the pentagon. We leave further details to the reader. 
 Proposition \ref{prop:aen-in-CM} is proved.

\paragraph{The category $\CMen$ and the PROB $\Ben$.}
Recall the PROB $\Ben$ governing graded bialgebras, see \S 
\ref{sec:main-res}\ref{par:uni-gr-bi}

\begin{thm}\label{thm:CM=B}
We have an equivalence of braided monoidal categories $\xi: \CMen\to\Ben$.
In particular, for any $n>0$ we have an equivalence of ordinary
(non-monoidal) categories $\xi_n: \CMen_n\to\Ben_n$. 
\end{thm}

\noindent{\sl Proof:} 
Recall that
 $\Ben$ has a graded bialgebra $\ba$. By definition, the components $\ba_n$
 are generating objects for $\Ben$, i.e., any other object is isomorphic to
 a tensor product of several of the $\ba_n$. 
 Similarly, the objects $[n]$ associated to $1\times 1$ matrices, are generating
 objects for $\CMen$. Indeed, any object $[M]$ associated to any $r\times s$ contingency matrix
 $M=\|m_{ij}\|$, is isomorphic to the tensor product  (in any order) of the individual $[m_{ij}]$,
the latter product being represented by the diagonal $rs\times rs$ matrix with entries $m_{ij}$
 in the corresponding order.  
 This  can be easily seen by moving the $m_{ij}$ around in the matrix by using anodyne 
 $\delta'$- and $\delta''$-isomoprhisms and their inverses.

 Next,  the graded bialgebra $\ba$, Corollary \ref{cor:bialg-CM-V-1}  and 
 Proposition \ref{prop:bialg-CM-V} give
  a braided monoidal functor
\[
\xi=\xi_\ba: \CMen\lra\Ben, \quad [M] \mapsto \ba_M.
\]
We prove that $\xi$ is an equivalence. For this, we use the graded
bialgebra $\aen$ in $\CMen$ constructed in Proposition \ref{prop:aen-in-CM}.
As $\ba\in \Ben$ is the universal graded bialgebra, we get a braided monoidal functor
\[
F=F_\aen: \Ben \lra\CMen, \quad \ba_n\mapsto \aen_n=[n]. 
\]
We claim that the functors $\xi$ and $F$ are quasi-inverse to each other.
Indeed, look at  the composition $F\xi: \CMen\to\CMen$, a braided
monoidal functor. It takes any generating object $[n]$  
to itself. Therefore $F\xi$ is isomorphic to $\Id$. Similarly, look at $\xi F: \Ben\to\Ben$.
It  is a braided monoidal functor which takes any generating object $\ba_n$ to itself.
Therefore $\xi F$ is isomorphic to $\Id$. \qed

\paragraph{Proof of Theorem  \ref{thm:main-sym}. }
Because of Theorem \ref{thm:CM=B}, Theorem  \ref{thm:main-sym}
can be
reformulated as follows.

\begin{refo}
For any abelian category $\Vc$ we have an equivalence of categories
\hfill\break 
$\Perv(\Sym^n(\CC);\Vc) \simeq \Fun(\CMen, \Vc)$. 
\end{refo}

This statement is a consequence (particular case) of the main result
of \cite {KS-hW} (Theorem 2.6) which describes perverse sheaves on
$W\backslash \hen$ where $\hen$ is the Cartan subalgebra of a reductive
complex Lie algebra $\gen$ and $W$ is the Weyl group of $\gen$. 
More precisely,  \cite {KS-hW}  deals with $\Vect_\k$-valued perverse sheaves,
but extension to perverse sheaves with values in an arbitrary abelian category $\Vc$ is trivial. 
Our case corresponds to $\gen=\gen\len_n$, when $\hen=\CC^n$ and
$W=S_n$, so $W\backslash \hen = \Sym^n(\CC)$. The description of
 \cite {KS-hW} is in terms of {\em mixed Bruhat sheaves} (Definition 2.1 there)
 which are certain diagrams with objects labelled by the set
 \[
 \Xi =\Xi_\gen \,=\, \bigsqcup_{I,J\subset\Delta_{\on{sim}}} W\backslash (W/W_I \times W/W_j)
 \]
where $\Delta_{\on{sim}}$ is the set of simple roots of $\gen$ and $W_I$, $I\subset
\Delta_{\on{sim}}$ is the subgroup in $W$ generated by the simple reflections
$s_\alpha$, $\alpha\in I$. For $\gen=\gen\len_n$ the set $\Xi_\gen$
is identified with $\CM_n$, see \cite {KS-cont}, and the axioms
of a mixed Bruhat sheaves become identical to the relations in the category
$\CMen$.  This finishes the proof.

\vfill\eject

\vskip 1cm

\small{

M.K.: Kavli IPMU, 5-1-5 Kashiwanoha, Kashiwa, Chiba, 277-8583 Japan. Email: 
\hfil\break
{\tt mikhail.kapranov@protonmail.com}

\smallskip

V.S.: Institut de Math\'ematiques de Toulouse, Universit\'e Paul Sabatier, 118 route de Narbonne, 
31062 Toulouse, France. Email: 
{\tt schechtman@math.ups-tlse.fr }

}

\ed